\documentclass[12pt]{amsart}
\usepackage{amssymb, upref}

\newcommand{\field}[1]{\mathbb{#1}}
\newcommand{\CC}{\field{C}}
\newcommand{\FF}{\field{F}}
\newcommand{\NN}{\field{N}}

\newcommand{\TT}{\field{T}}
\newcommand{\ZZ}{\field{Z}}

\newcommand{\Ff}{\mathcal F}

\newcommand{\Ii}{\mathcal I}
\newcommand{\Kk}{\mathcal K}
\newcommand{\Ll}{\mathcal L}
\newcommand{\Mm}{\mathcal M}
\newcommand{\Oo}{\mathcal O}
\newcommand{\Tt}{\mathcal T}

\newcommand{\MCE}{\operatorname{MCE}}

\newcommand{\id}{\operatorname{id}}
\newcommand{\lsp}{\operatorname{span}}
\newcommand{\clsp}{\operatorname{\overline{span\!}\,\,}}

\newcommand{\aF}{\widetilde{F}}
\newcommand{\Augl}{\widetilde{l}}
\newcommand{\Xle}[1]{\widetilde{X}_{#1}}
\newcommand{\iotale}{\widetilde{\iota}\hskip1pt}
\newcommand{\phile}{\widetilde{\phi}}
\newcommand{\Gopp}{\Gamma^{\operatorname{opp}}}
\newcommand{\intfrm}[1]{\Pi{#1}}

\newcommand{\M}{M}
\newcommand{\Tc}{\Tt_{\rm cov}}
\newcommand{\NO}[1]{\operatorname{\mathcal{NO}}_{\!#1}}
\newcommand{\fixNO}[2]{\NO{#1}^{\hskip0.1em#2}}
\newcommand{\Sbit}{S}
\newcommand{\notSbit}{P \setminus S}

\theoremstyle{plain}
\newtheorem{theorem}{Theorem}[section]
\newtheorem*{theorem*}{Theorem}
\newtheorem*{prop*}{Proposition}
\newtheorem{cor}[theorem]{Corollary}
\newtheorem{lemma}[theorem]{Lemma}
\newtheorem{prop}[theorem]{Proposition}

\theoremstyle{remark}
\newtheorem{rmk}[theorem]{Remark}

\newtheorem{example}[theorem]{Example}

\theoremstyle{definition}
\newtheorem{dfn}[theorem]{Definition}

\newtheorem{notation}[theorem]{Notation}

\numberwithin{equation}{section}

\title[$C^*$-algebras associated to product systems]
{$C^*$-algebras associated to product
systems of Hilbert bimodules}
\author{Aidan Sims}
\address{Aidan Sims\\
School of Mathematics and Applied Statistics\\
University of Wollongong\\
NSW 2522\\
AUSTRALIA} \email{aidan\_sims@uow.edu.au}

\author{Trent Yeend}
\address{Trent Yeend\\
School of Mathematical and Physical Sciences\\
University of Newcastle\\
NSW  2308\\
AUSTRALIA} \email{trent.yeend@newcastle.edu.au}
\date{December 19, 2007; revised January 6, 2009.}
\thanks{This research was supported by the Australian Research Council.}

\begin{document}

\begin{abstract}
Let $(G,P)$ be a quasi-lattice ordered group and let $X$ be a
compactly aligned product system over $P$ of Hilbert bimodules
in the sense of Fowler. Under mild hypotheses we associate to
$X$ a $C^*$-algebra which we call the Cuntz-Nica-Pimsner
algebra of $X$. Our construction generalises a number of
others: a sub-class of Fowler's Cuntz-Pimsner algebras for
product systems of Hilbert bimodules; Katsura's formulation of
Cuntz-Pimsner algebras of Hilbert bimodules; the $C^*$-algebras
of finitely aligned higher-rank graphs; and Crisp and Laca's
boundary quotients of Toeplitz algebras. We show that for a
large class of product systems $X$, the universal
representation of $X$ in its Cuntz-Nica-Pimsner algebra is
isometric.
\end{abstract}

\maketitle

\section{Introduction}

In this article we introduce and begin to analyse a class of
$C^*$-algebras, which we call Cuntz-Nica-Pimsner algebras,
associated to product systems of Hilbert bimodules. This work
draws on and generalises a substantial body of previous work in
a number of related areas: results of \cite{FMR, FR, K6, MS,
Pim97} on $C^*$-algebras associated to Hilbert bimodules; the
study of $C^*$-algebras associated to product systems in
\cite{Din, Din2, F99,FR98}; the theory of $C^*$-algebras
associated to higher-rank graphs \cite{KP,RS1,RSY2}; and the
theory of Toeplitz algebras (and quotients thereof) associated
to quasi-lattice ordered groups \cite{CL1, CL, LR, Nica}.
Consequently, putting our results in context and indicating
their significance requires some discussion.

\subsection{$C^*$-algebras associated to Hilbert bimodules}

In \cite{Pim97}, Pimsner associated to each Hilbert $A$--$A$
bimodule $X$ two $C^*$-algebras $\Oo_X$ and $\Tt_X$. He showed
that, as the notation suggests, the $C^*$-algebras $\Oo_X$
generalise the Cuntz-Krieger algebras $\Oo_A$ associated to
$\{0,1\}$-matrices $A$ in \cite{CK}, and the algebras $\Tt_X$
generalise their Toeplitz extensions $\Tt\Oo_A$. According to
\cite{Pim97}, a representation of $X$ in a $C^*$-algebra $B$ is
a pair $(\pi,\psi)$ where $\pi : A \to B$ is a
$C^*$-homomorphism, $\psi : X \to B$ is linear, and the pair
carries the Hilbert $A$--$A$ bimodule structure on $X$ to the
natural Hilbert $B$--$B$ bimodule structure on $B$. Pimsner
proved that $\Tt_X$ is universal for representations of $X$,
and that $\Oo_X$ is universal for representations satisfying a
covariance condition now known as Cuntz-Pimsner covariance, and
in particular is a quotient of $\Tt_X$.

In the spirit of Coburn's Theorem for the classical Toeplitz
algebra \cite{Coburn:isometryI}, Fowler and Raeburn proved a
uniqueness theorem for $\Tt_X$ \cite[Theorem~3.1]{FR}. A key
example in their work, and an informative application of
Pimsner's ideas, is related to graph $C^*$-algebras. A directed
graph consists of a set $E^0$ of vertices, a set $E^1$ of
edges, and maps $r,s : E^1 \to E^0$, called the range and
source maps, which give the edges their direction: the edge $e$
is directed from $s(e)$ to $r(e)$. When we discuss graphs and
their $C^*$-algebras here, we follow the conventions of
\cite{CBMSbk}.

Fowler and Raeburn used their analysis of $\Tt_X$ to
investigate graph algebras by associating a Hilbert
$c_0(E^0)$--$c_0(E^0)$ bimodule $X(E)$ to each directed graph
$E$. Previously \cite{KPR, KPRR} $C^*$-algebras had been
associated to directed graphs $E$ which are row-finite (each
vertex receives at most finitely many vertices) and have no
sources (each vertex receives at least one edge). In~\cite{FR},
Fowler and Raeburn showed that for such graphs, $\Tt_{X(E)}$
can naturally be thought of as a Toeplitz extension of the
graph $C^*$-algebra $C^*(E)$ of \cite{KPR}. They observed that
$E$ is row-finite precisely when the homomorphism $\phi : A \to
\Ll(X(E))$ which implements the left action takes values in the
algebra $\Kk(X(E))$ of generalised compact operators on $X(E)$
and that $E$ has no sources precisely when $\phi$ is injective.
As Pimsner's theory does not require that the left action be by
compact operators, Fowler and Raeburn's results suggested what
is now the accepted definition of the graph $C^*$-algebra of a
non-row-finite graph \cite{FLR}. Results of Exel and Laca
\cite{EL} were used to prove a version of the Cuntz-Krieger
uniqueness theorem for arbitrary graph $C^*$-algebras in
\cite{FLR}, and direct methods were used to extend a number of
other graph $C^*$-algebraic results to the non-row-finite
setting in \cite{BHRS}.

Via the connection between graph $C^*$-algebras and
Cuntz-Pimsner algebras discussed above, the uniqueness theorems
of \cite{BHRS, FLR} suggested an alternate approach to
Cuntz-Pimsner algebras when $\phi$ is not injective.
Specifically, when interpreted in terms of representations of
the bimodule $X(E)$, the gauge-invariant uniqueness theorem of
\cite{BHRS} suggested that one could weaken Pimsner's
covariance condition to obtain a covariance condition for which
the universal $C^*$-algebra satisfies the following two
criteria: that the universal representation is injective on $A$
(in Pimsner's theory, this requires that $\phi$ is injective);
and that any representation of the universal $C^*$-algebra
which respects the gauge action of $\TT$ and is injective on
$A$ is faithful.

In \cite{K6}, Katsura identified such a covariance condition: it is
the defining relation for a relative Cuntz-Pimsner algebra (see
\cite{FMR, MS}) with respect to a certain ideal of $A$. Katsura's
universal algebra $\Oo_X$ satisfies the two criteria set forth in the
preceding paragraph, and under Katsura's definition, $\Oo_{X(E)}
\cong C^*(E)$ for arbitrary graphs $E$.

\subsection{$C^*$-algebras associated to product
systems}\label{sec:intro ps}

Let $P$ be a semigroup with identity $e$. Informally, a product
system over $P$ of Hilbert $A$--$A$ bimodules is a semigroup $X
= \bigsqcup_{p \in P} X_p$ such that each $X_p$ is a
right-Hilbert $A$--$A$ bimodule, and $x \otimes_A y \mapsto xy$
determines an isomorphism of $X_p \otimes_A X_q$ onto $X_{pq}$
for all $p,q \in P\setminus\{e\}$. These objects were
introduced in this generality by Fowler in \cite{F99} as
generalisations of the continuous product systems of Hilbert
spaces introduced by Arveson in \cite{Arv} and their discrete
analogues introduced by Dinh in \cite{Din, Din2}.

In \cite{F99}, Fowler considered a class of product systems $X$
over semigroups $P$ arising in quasi-lattice ordered groups
$(G,P)$ which he calls \emph{compactly aligned} product systems
(see~\eqref{eq:compactly aligned}). Inspired by work of Nica
\cite{Nica} and of Laca and Raeburn \cite{LR} on Toeplitz
algebras associated to quasi-lattice ordered groups, Fowler
introduced and studied what he called Nica covariant
representations of $X$ and the associated universal
$C^*$-algebra $\Tc(X)$. When $P = \NN$ Nica covariance is
automatic and $\Tc(X)$ coincides with Pimsner's $\Tt_{X_1}$.
Fowler's main theorem \cite[Theorem~7.2]{F99} gave a spatial
criterion for faithfulness of a representation of $\Tc(X)$.
This theorem generalised both Laca and Raeburn's uniqueness
theorem~\cite[Theorem~3.7]{LR} for $C^*(G,P)$ and Fowler and
Raeburn's uniqueness theorem~\cite[Theorem~3.1]{FR} for the
Toeplitz algebra of a single Hilbert bimodule.

In~\cite{F99} Fowler also proposed a Cuntz-Pimsner covariance
condition for representations of compactly aligned product
systems and an associated universal $C^*$-algebra $\Oo_X$.
Fowler defined a representation $\psi$ of a product system $X$
to be Cuntz-Pimsner covariant if its restriction to each $X_p$
is Cuntz-Pimsner covariant in Pimsner's sense. Fowler showed
that his Cuntz-Pimsner covariance condition implies Nica
covariance under the hypotheses that: each pair of elements of
$P$ has a common upper bound; each $X_p =
\overline{\phi_p(A)X_p}$; and each $\phi_p(A) \subset
\Kk(X_p)$.

Results of \cite{RS1} generalised the construction of a
bimodule $X(E)$ from a graph $E$ to the construction of a
product system of bimodules $X(\Lambda)$ over $\NN^k$ from a
$k$-graph $\Lambda$. In particular, \cite[Theorem~5.4]{RS1}
identifies the $k$-graphs $\Lambda$ for which $X(\Lambda)$ is
compactly aligned; such $k$-graphs are said to be finitely
aligned. An analysis of the Toeplitz algebra $\Tt C^*(\Lambda)$
of such a finitely-aligned $k$-graph $\Lambda$ based on
Fowler's results \cite{RS1} led to the formulation in
\cite{RSY2} of a Cuntz-Krieger relation for finitely aligned
$k$-graphs. Direct methods were used in \cite{RSY2} to prove
versions of the standard uniqueness theorems for higher-rank
graph $C^*$-algebras.

Just as the uniqueness theorems of~\cite{BHRS} informed the
work of Fowler, Muhly and Raeburn~\cite{FMR} and of Katsura
\cite{K5,K6}, the uniqueness theorems of \cite{RSY2} provide an
informative model for $C^*$-algebras associated to compactly
aligned product systems. In particular, the gauge-invariant
uniqueness theorem of \cite{RSY2} prompts us to seek a
$C^*$-algebra $\NO{X}$ satisfying two criteria:
\begin{itemize}
\item[(A)] the universal homomorphism $j_X$ restricts to an
    injection on $X_e = A$; and
\item[(B)] any representation of $\NO{X}$ which is faithful
    on $j_X(A)$ is faithful on the fixed-point algebra
    $\fixNO{X}{\delta}$ (where $\delta$ is the canonical
    gauge coaction of the enveloping group $G$ of $P$ on
    $\NO{X}$).
\end{itemize}
To see the analogy of criterion~(B) with the corresponding
criterion given above for the Cuntz-Pimsner algebra of a single
Hilbert module, consider the case where $P = \NN^k$. Averaging
over the gauge action $\delta$ of $\TT^k$ on $\NO{X}$
determines a faithful conditional expectation of $\NO{X}$ onto
$\fixNO{X}{\delta}$. A standard argument then proves that every
representation which is faithful on $\fixNO{X}{\delta}$ and
respects $\delta$ is faithful on all of $\NO{X}$. So our
criterion~(B) implies that if a representation of $\NO{X}$
restricts to an injection of $j_X(A)$ and respects $\delta$,
then it is injective.

Kumjian and Pask's results regarding row-finite $k$-graphs with
no sources suggest that we can expect Fowler's $\Oo_X$ to
satisfy (A)~and~(B) when the left action on each fibre is
injective and by compact operators and $P = \NN^k$. However,
$\Oo_X$ will not always fit the bill. If $a \in A$ acts
trivially on the left of some $X_p$, then Fowler's relation
forces $j_X(a) = 0$. Moreover examples of
\cite[Appendix~A]{RSY2}, when interpreted in terms of product
systems, show that if $A$ does not act compactly on the left of
every $X_p$, Fowler's $\Oo_X$ may not satisfy criterion~(B).
Indeed, the examples of \cite[Appendix~A]{RSY2} show that the
same problems persist even if, in Fowler's definition of
$\Oo_X$, we replace Pimsner's covariance condition with
Katsura's. The situation is less clear when pairs in $P$ need
not have a common upper bound because in this case we have
fewer guiding examples. We take the approach, different from
Fowler's, that $\NO{X}$ should always be a quotient of
$\Tc(X)$. The situation is then clearer: the $C^*$-algebra
which is universal for representations satisfying both Nica
covariance and Fowler's Cuntz-Pimsner covariance will not
satisfy~(A) if there exist $p,q \in P$ with $p \vee q =
\infty$.

The approach that $\NO{X}$ should be a quotient of $\Tc(X)$ is
justified by examples in isometric representation theory for
quasi-lattice ordered groups. Fowler~\cite{F99} showed that
when each $X_p$ is a $1$-dimensional Hilbert space and
multiplication in $X$ is implemented by multiplication of
complex numbers, his $\Tc(X)$ agrees with the $C^*$-algebra
$C^*(G,P)$ universal for Nica covariant representations of $P$
\cite{CL1, LR, Nica}. Crisp and Laca have recently studied what
they call boundary quotients of Toeplitz algebras associated to
quasi-lattice ordered groups \cite{CL}. Their results show that
the relationship between their boundary quotient and $C^*(G,P)$
is often analogous to the relationship between the Cuntz
algebras and their Toeplitz extensions. In particular if $G =
\FF_n$ is the free group on $n$ generators, and $P = \FF^+_n$
is its positive cone, then $C^*(G,P)$ is isomorphic to the
Toeplitz extension $\Tt\Oo_n$ of $\Oo_n$ and Crisp and Laca's
boundary quotient is isomorphic to $\Oo_n$ itself.

More generally, \cite[Theorem~6.7]{CL} shows that when $G$ is a
right-angled Artin group with trivial centre (and $P$ is its
positive cone), the associated boundary quotient is simple and
purely infinite. Regarded as a statement about a $C^*$-algebra
associated to the product system over $P$ with $1$-dimensional
fibres, simplicity is equivalent to~(B). Hence boundary
quotients provide important motivation and a good test-case for
our theory when $G \not= \ZZ^k$. In particular, the results of
Section~\ref{sec:CrispLaca} suggest that we are on the right
track for fairly general quasi-lattice ordered groups.

\subsection{Outline of the paper.}

In this paper we combine ideas of \cite{F99, FMR, K6} with
intuition drawn from the theory of $k$-graph $C^*$-algebras
\cite{RSY2} to associate what we call a Cuntz-Nica-Pimsner
algebra $\NO{X}$ to a broad class of compactly aligned product
systems $X$ of Hilbert bimodules. This $\NO{X}$ is universal
for a class of representations which we refer to as
Cuntz-Nica-Pimsner covariant (or CNP-covariant for short). By
definition these are the representations which are Nica
covariant in the sense of Fowler, and also satisfy a
Cuntz-Pimsner covariance relation which looks substantially
different from Fowler's (but agrees with Fowler's under a
number of additional hypotheses on $X$). Our ultimate aim is to
verify that $\NO{X}$ satisfies both (A)~and~(B) above.

Our main result, Theorem~\ref{thm:isometric}, shows that our
$C^*$-algebra $\NO{X}$ satisfies~(A). In Sections
\ref{sec:Katsura's}, \ref{sec:hrg}~and~\ref{sec:CrispLaca}, we
present evidence that it also satisfies~(B). In
Section~\ref{sec:Katsura's}, we prove that our construction
agrees with Katsura's when $P = \NN$, and in
Section~\ref{sec:hrg}, we show that given a $k$-graph $\Lambda$
and the associated product system $X(\Lambda)$ as constructed
in \cite{RS1}, our $\NO{X}$ coincides with the Cuntz-Krieger
algebra $C^*(\Lambda)$ of \cite{RSY2}. The gauge-invariant
uniqueness theorems
\cite[Theorem~4.2]{RSY2}~and~\cite[Theorem~6.2]{K6} then imply
that our definition satisfies~(B). In
Section~\ref{sec:CrispLaca} we prove that our Cuntz-Pimsner
covariance relation is compatible with the defining relations
for Crisp and Laca's boundary quotients~\cite{CL}. In
particular, Crisp and Laca's uniqueness theorem
\cite[Theorem~6.7]{CL} shows that $\NO{X}$ satisfies~(B) when
$G$ is a right-angled Artin group with trivial centre and $X$
is the product system over $P$ with $1$-dimensional Hilbert
spaces for fibres and multiplication implemented by
multiplication of complex numbers.

We show in Section~\ref{sec:Fowler's} that $\NO{X}$ coincides with
Fowler's $\Oo_X$ when \cite[Proposition~5.4]{F99} suggests that it
might --- namely when each pair in $P$ has a common upper bound, and
each $\phi_p$ is injective with $\phi_p(A) \subset \Kk(X_p)$.

\subsection*{Acknowledgements}

It came to our attention in the late stages of preparation of
this manuscript that Toke Carlsen, Nadia Larsen and Sean
Vittadello were also considering the question of a satisfactory
notion of Cuntz-Pimsner covariance for representations of
product systems. We thank them for many helpful and insightful
comments on this manuscript which have significantly clarified
both our ideas and our exposition of them.

We also thank Iain Raeburn for a number of helpful discussions.

\section{Definitions}

In this section, we recall the definitions and notation
described in \cite[Section~1]{F99}.

\subsection{Hilbert bimodules}

We attempt to summarise only those aspects of Hilbert bimodules of
direct relevance to this paper. We refer the reader to \cite{Black,
Lan, TFB} for more detail.

Let $A$ be a $C^*$-algebra, and let $X$ be a complex vector
space carrying a right action of $A$. Suppose that $\langle
\cdot, \cdot \rangle_A : X \times X \to A$ is linear in the
second variable and conjugate linear in the first variable and,
for $x,y \in X$ and $a \in A$, satisfies:
\begin{enumerate}
\item $\langle x, y\rangle_A = \langle y,x\rangle_A^*$;
\item $\langle x, y \cdot a\rangle_A = \langle x,y\rangle_A
    a$;
\item $\langle x, x \rangle_A$ is a positive element of
    $A$; and
\item $\langle x, x \rangle_A = 0 \iff x = 0$.
\end{enumerate}
The formula $\|x\| = \|\langle x, x\rangle_A\|^{\frac12}$ defines a
norm on $X$. If $X$ is complete in this norm, we call it a
\emph{right-Hilbert $A$-module}.

Let $X_A$ be a right-Hilbert $A$-module. A map $T : X \to X$ is
said to be \emph{adjointable} if there is a map $T^* : X \to X$
such that $\langle Tx, y\rangle_A = \langle x, T^*y\rangle_A$
for all $x,y \in X$. Every adjointable operator on $X$ is
norm-bounded and linear, and the adjoint $T^*$ is unique. The
collection $\Ll(X)$ of adjointable operators on $X$ endowed
with the operator norm is a $C^*$-algebra. Given $x,y \in X$,
there is an adjointable operator $x \otimes y^*$ on $X$
determined by the formula $(x \otimes y^*)(z) = x \cdot \langle
y,z \rangle_A$. We call operators of this form
\emph{generalised rank-1 operators}. The subspace $\Kk(X) :=
\clsp\{x \otimes y^* : x,y \in X\}$ is an essential ideal of
$\Ll(X)$ whose elements we refer to as \emph{generalised
compact operators on $X$}.

A \emph{right-Hilbert $A$--$A$ bimodule} is a right-Hilbert $A$
module together with a homomorphism $\phi : A \to \Ll(X)$. We
think of $\phi$ as implementing a left action of $A$ on $X$, so
we typically write $a \cdot x$ for $\phi(a)x$. Because $\phi(a)
\in \Ll(X)$ for all $a \in A$, we automatically have $a \cdot
(x \cdot b) = (a \cdot x) \cdot b$ for all $a,b \in A$ and $x
\in X$.

An important special case is the bimodule $_AA_A$ with inner
product given by $\langle a,b\rangle_A = a^*b$ and right- and
left-actions given by multiplication in $A$. The $C^*$-algebra
$\Ll(_AA_A)$ is isomorphic to the multiplier algebra $\Mm(A)$,
and the homomorphism that takes $a \in A$ to
left-multiplication by $a$ on $X$ is an isomorphism of $A$ onto
$\Kk(_AA_A)$.

We form the balanced tensor product $X \otimes_A Y$ of two
right-Hilbert $A$--$A$ bimodules as follows. Let $X \odot Y$ be the
algebraic tensor product of $X$ and $Y$ as complex vector spaces. Let
$X \odot_A Y$ be the quotient of $X \odot Y$ by the subspace spanned
by vectors of the form $x\cdot a \odot y - x \odot a\cdot y$ where $x
\in X$, $y \in Y$ and $a \in A$. The formula
\[
\langle x_1 \odot y_1, x_2 \odot y_2 \rangle_A = \langle y_1,
\langle x_1, x_2\rangle_A \cdot y_2\rangle_A
\]
determines a bounded sesquilinear form on $X \odot_A Y$. Let $N
= \lsp\{n \in X \odot_A Y : \langle n, n\rangle_A = 0_A\}$.
Then $\|z + N\| = \inf\{\|\langle z + n, z + n\rangle_A\|^{1/2}
: n \in N\}$ defines a norm on $(X \odot_A Y)/N$, and $X
\otimes_A Y$ is the completion of $(X \odot_A Y)/N$ in this
norm.

If $X$ and $Y$ are right-Hilbert $A$--$A$ bimodules and $S \in
\Ll(X)$, then there is an adjointable operator $S \otimes 1_Y$
(with adjoint $S^* \otimes 1_Y$) on $X \otimes_A Y$ determined
by $(S \otimes 1_Y)(x \otimes_A y) = Sx \otimes_A y$ for all $x
\in X$ and $y \in Y$. The formula $a \mapsto \phi(a) \otimes
1_Y$ therefore determines a homomorphism of $A$ into $\Ll(X
\otimes_A Y)$. The notation $X \otimes_A Y$ always refers to
the right-Hilbert $A$--$A$ bimodule in which the left action is
implemented by this homomorphism.

\subsection{Semigroups and product systems of Hilbert bimodules}

Let $P$ be a discrete multiplicative semigroup with identity
$e$, and let $A$ be a $C^*$-algebra. A \emph{product system
over $P$ of right-Hilbert $A$--$A$ bimodules} is a semigroup $X
= \bigsqcup_{p \in P} X_p$ such that:
\begin{itemize}
\item[(1)] for each $p \in P$, $X_p \subset X$ is a right-Hilbert
    $A$--$A$ bimodule;
\item[(2)] the identity fibre $X_e$ is equal to the bimodule
    $_AA_A$;
\item[(3)] for $p,q \in P\setminus\{e\}$ there is an isomorphism
    $\M_{p,q} : X_p \otimes_A X_q \to X_{pq}$ satisfying
    $M_{p,q}(x \otimes_A y) = xy$ for all $x \in X_p$ and $y \in
    X_q$; and
\item[(4)] multiplication in $X$ by elements of $X_e = A$
    implements the actions of $A$ on each $X_p$; that is
    $ax = a \cdot x$ and $xa = x \cdot a$ for all $p \in
    P$, $x \in X_p$ and $a \in X_e$.
\end{itemize}
For $p \in P$, we denote the homomorphism of $A$ to $\Ll(X_p)$
which implements the left action by $\phi_p$, and we denote the
$A$-valued inner product on $X_p$ by $\langle \cdot, \cdot
\rangle_A^p$.

By (2)~and~(4), for $p \in P$, multiplication in $X$ induces
maps $\M_{p,e} : X_p \otimes_A X_e \to X_p$ and $M_{e,p} : X_e
\otimes_A X_p \to X_p$ as in~(3). Each $M_{p,e}$ is
automatically an isomorphism by~\cite[Corollary~2.7]{TFB}. We
do not insist that $M_{e,p}$ is an isomorphism as this is too
restrictive if we want to capture Pimsner's theory (which does
not require that $\overline{\phi(A)X} = X$).

Because multiplication in $X$ is associative, we have
$\phi_{pq}(a)(xy) = (\phi_p(a)x)y$ for all $x \in X_p$, $y \in
X_q$ and $a \in A$.

Given $p, q \in P$ with $p \not= e$, the isomorphism $M_{p,q} :
X_p \otimes_A X_q \to X_{pq}$ allows us to define a
homomorphism $\iota^{pq}_p : \Ll(X_p) \to \Ll(X_{pq})$ by
\[
\iota^{pq}_p(S) = \M_{p, q} \circ (S \otimes 1_{X_{q}})
\circ \M_{p, q}^{-1}.
\]
We may alternatively characterise $\iota^{pq}_p$ by the formula
$\iota^{pq}_p(S)(xy) = (Sx)y$ for all $x \in X_p$, $y \in X_{q}$ and
$S \in \Ll(X_p)$. When $p = e$, we do not have $X_p \otimes X_q \cong
X_{pq}$; however, since $X_p = A$, we may define $\iota^q_e$ on
$\Kk(X_e) \cong A$ by $\iota^q_e(a) = \phi_q(a)$ for all $a \in A$.
As a notational convenience, if $p, r \in P$ and $r \not= pq$ for any
$q \in P$, we define $\iota^r_p : \Ll(X_p) \to \Ll(X_r)$ to be the
zero map $\iota^r_p(S) = 0_{\Ll(X_r)}$ for all $S \in \Ll(X_p)$.

We will primarily be interested in semigroups $P$ of the
following form. Following Nica \cite{Nica}, we say that $(G,P)$
is a \emph{quasi-lattice ordered group} if: $G$ is a discrete
group and $P$ is a subsemigroup of $G$; $P \cap P^{-1} =
\{e\}$; and with respect to the partial order $p \le q \iff
p^{-1}q \in P$, any two elements $p,q \in G$ which have a
common upper bound in $P$ have a least upper bound $p \vee q
\in P$. We write $p \vee q = \infty$ to indicate that $p,q \in
G$ have no common upper bound in $P$, and we write $p \vee q <
\infty$ otherwise.

Let $(G,P)$ be a quasi-lattice ordered group, and let $X$ be a
product system over $P$ of right-Hilbert $A$--$A$ bimodules. We
say that $X$ is \emph{compactly aligned} if
\begin{equation}\label{eq:compactly aligned}
\parbox{0.9\textwidth}{for all $p,q \in P$
such that $p \vee q < \infty$, and for all $S \in \Kk(X_p)$ and
$T \in \Kk(X_q)$, we have $\iota^{p \vee q}_p(S) \iota^{p \vee
q}_q(T) \in \Kk(X_{p\vee q})$.}
\end{equation}
Note that this condition does not imply compactness of either
$\iota^{p \vee q}_p(S)$ or $\iota^{p \vee q}_q(T)$.

\subsection{Representations of product systems and Nica covariance}
Let $(G,P)$ be a quasi-lattice ordered group, and let $X$ be a
compactly aligned product system over $P$ of right-Hilbert
$A$--$A$ bimodules. Let $B$ be a $C^*$-algebra, and let $\psi$
be a function from $X$ to $B$. For $p \in P$, let $\psi_p =
\psi|_{X_p}$. We call $\psi$ a \emph{representation of $X$} if
\begin{itemize}
\item[(${\rm T1}$)] each $\psi_p : X_p \to B$ is linear,
    and $\psi_e$ is a $C^*$-homomorphism;
\item[(${\rm T2}$)] $\psi_p(x)\psi_q(y) = \psi_{pq}(xy)$ for all
    $p,q \in P$, $x \in X_p$, and $y \in X_q$; and
\item[(${\rm T3}$)] $\psi_e(\langle x, y\rangle^p_A) =
    \psi_p(x)^* \psi_p(y)$ for all $p \in P$, and $x,y \in X_p$.
\end{itemize}

\begin{rmk}
Our definition agrees with Fowler's \cite[Definition~2.5]{F99}:
condition~(4) of the definition of a product system together
with (T1)--(T3) ensures that each $(\psi_e, \psi_p)$ is
representation of $X_p$ in the sense of Pimsner. It then
follows from Pimsner's results (see \cite[p.~202]{Pim97}) that
for each $p \in P$ there is a homomorphism $\psi^{(p)} :
\Kk(X_p) \to B$ which satisfies $\psi^{(p)}(x \otimes y^*) =
\psi_p(x)\psi_p(y)^*$ for all $x,y \in X_p$.
\end{rmk}

We say that a representation $\psi$ of $X$ is \emph{Nica
covariant} if
\begin{itemize}
\item[(${\rm N}$)] For all $p,q \in P$ and all $S \in
    \Kk(X_p)$ and $T \in \Kk(X_q)$,
\[
\displaystyle \psi^{(p)}(S)\psi^{(q)}(T) =
\begin{cases}
\psi^{(p\vee q)}\big(\iota^{p\vee q}_p(S)\iota^{p\vee q}_q(T)\big)
& \text{if $p\vee q < \infty$} \\
0 &\text{otherwise.}
\end{cases}
\]
\end{itemize}

Results of \cite{F99} show that there exist a $C^*$-algebra
$\Tc(X)$ and a Nica covariant representation $i$ of $X$ in
$\Tc(X)$ which are universal in the sense that:
\begin{itemize}
\item[(1)] $\Tc(X)$ is generated by $\{i(x) : x \in X\}$;
    and
\item[(2)] if $\psi$ is any Nica covariant representation
    of $X$ on a $C^*$-algebra $B$ then there is a unique
    homomorphism $\psi_* : \Tc(X) \to B$ such that $\psi_*
    \circ i = \psi$.
\end{itemize}

\section{Cuntz-Nica-Pimsner covariance and the Cuntz-Nica-Pimsner algebra}
In this section we present our definition of Cuntz-Pimsner
covariance for a compactly aligned product system $X$. We then
introduce the Cuntz-Nica-Pimsner algebra $\NO{X}$.

To present our definition we begin by introducing a collection
of bimodules $\Xle{p}$, $p \in P$ which we associate to a
product system $X$ over $P$. The $\Xle{p}$ do not form a
product system; rather, they play a r\^{o}le similar to that
played by the sets $\Lambda^{\le m}$ for a locally convex
$k$-graph in \cite{RSY1}. We first recall some standard
notation for direct sums of Hilbert modules.

Let $J$ be a set, let $X_j$ be a right-Hilbert $A$--$A$
bimodule for each $j \in J$, and let $X = \bigsqcup_{j \in J}
X_j$. Let $\Gamma_c(J, X)$ be the space of all
finitely-supported sections $x : J \to X$; that is $x(j) \in
X_j$ for all $j \in J$. The formula $\langle x, y \rangle_A =
\sum_{j \in J} \langle x(j), y(j) \rangle_A$ defines an
$A$-valued inner-product on $\Gamma_c(J, X)$. The completion of
$\Gamma_c(J, X)$ in the norm arising from this inner product is
called the \emph{direct sum} of the $X_j$ and denoted
$\bigoplus_{j \in J} X_j$. Endowed with pointwise left- and
right-actions of $A$, $\bigoplus_{j \in J} X_j$ is itself a
right-Hilbert $A$--$A$ bimodule.

Let $X$ be a right-Hilbert $A$ module, and let $I$ be an ideal
of $A$. Let $X \cdot I$ denote $\{x \cdot a : x \in X, a \in
I\}$. It is well-known that
\begin{equation}\label{eq:Hewitt-Cohen}
X \cdot I = \{x \in X : \langle x, x \rangle_A \in I\} = \clsp\{x \cdot a : x \in X, a \in I\}.
\end{equation}
One way to see this is as follows. First, if $x \in X$ and $a
\in I$, then $\langle x \cdot a, x \cdot a \rangle_A =
a^*(\langle x, x \rangle_A) a \in I$. Now suppose $\langle x, x
\rangle_A \in I$. Then the element $a = (\langle x, x
\rangle_A)^{1/4}$ also belongs to $I$. Since
\cite[Lemma~4.4]{Lan} implies that $x = w \cdot a$ for some $w
\in X$, we then have $x \in X \cdot I$. We have now established
the first equality in~\eqref{eq:Hewitt-Cohen}. The second
follows from \cite[Lemma~3.23]{TFB}: though this lemma is
stated for imprimitivity bimodules, the proof of the assertion
we are using requires only that $X$ is a right-Hilbert
$A$-module.

Since $I$ is an ideal, \eqref{eq:Hewitt-Cohen} implies that $X
\cdot I$ is itself a right-Hilbert $A$--$A$ bimodule.

\begin{dfn}\label{dfn:Xle def}
Let $(G,P)$ be a quasi-lattice ordered group, and let $X$ be a
product system over $P$ of right-Hilbert $A$--$A$ bimodules.
Define $I_e = A$, and for $p \in P$ define $I_p = \bigcap_{e <
r \le p} \ker(\phi_r) \lhd A$. For $q \in P$, we define the
right-Hilbert $A$--$A$ bimodule $\Xle{q}$ by
\[
\Xle{q} = \bigoplus_{p \le q} X_p \cdot
I_{p^{-1}q}.
\]
We write $\phile_q$ for the homomorphism from $A$ to
$\Ll(\Xle{q})$ which implements the left action of $A$ on
$\Xle{q}$. That is, $(\phile_q(a)x)(p) = \phi_p(a)x(p)$ for $p
\le q$.
\end{dfn}

\begin{lemma}\label{lem:Xle characterisation}
Let $(G,P)$ be a quasi-lattice ordered group, and let $X$ be a
product system over $P$ of right-Hilbert $A$--$A$ bimodules. Fix $p,
q \in P$ with $p \le q$. For $x \in X_p$, we have
\begin{equation}\label{eq:equiv conditions}
x \in X_p \cdot I_{p^{-1}q}
 \qquad\iff\qquad
xy = 0 \text{ for all $e < r \le p^{-1}q$ and $y \in X_r$.}
\end{equation}
\end{lemma}
\begin{proof}
If $p = e$, then $X_p = A$, and $xy = \phi_r(x)y$ for $r \le
e^{-1}q = q$ and $y \in X_r$. Hence both sides
of~\eqref{eq:equiv conditions} reduce to $x \in I_q$.

Suppose $p \not= e$. By~\eqref{eq:Hewitt-Cohen}
\begin{align}
x \in X_p \cdot I_{p^{-1}q}
 &\iff  \langle x, x \rangle^p_A \in I_{p^{-1}q} \nonumber \\
 &\iff \langle x, x \rangle^p_A \in \ker(\phi_r)\text{ for all $e < r \le p^{-1}q$.} \label{eq:in ker(phi_r)}
\end{align}
Remark~2.29 of \cite{TFB} implies that a positive adjointable
operator $T$ on a right-Hilbert $A$-module $X$ is equal to zero if
and only if $\langle z, Tz \rangle_A = 0$ for all $z \in X$. Since
$\langle x, x \rangle^p_A$ is a positive element of $A$, its image
under $\phi_r$ is positive for any $r$, so~\eqref{eq:in ker(phi_r)}
implies that
\begin{equation}\label{eq:ip characterisation}
x \in X_p \cdot I_{p^{-1}q}\quad
 \iff\quad \langle y, \langle x, x \rangle^p_A \cdot y \rangle^r_A = 0\text{ for all $e < r \le p^{-1}q$ and $y \in X_r$}.
\end{equation}
By definition of the inner product on the internal tensor product of
Hilbert bimodules, for $e < r \le p^{-1}q$ and $y \in X_r$, we have
\[
\langle y, \langle x, x \rangle^p_A \cdot y \rangle^r_A
 = \langle x \otimes_A y, x \otimes_A y \rangle_A
 = \|x \otimes_A y\|^2.
\]
Combining this with~\eqref{eq:ip characterisation}, we have
\[
x \in X_p \cdot I_{p^{-1}q}
\qquad\iff\qquad
 x \otimes_A y = 0\text{ for all $e < r \le p^{-1}q$ and all $y \in X_r$.}
\]
Since $p \not= e$, each $\M_{p,r}$ is an isomorphism, and the
result follows.
\end{proof}

\begin{example}\label{ex:Xle<->Lambdale}
Let $\Lambda$ be a finitely aligned $k$-graph and let $X =
X(\Lambda)$ be the corresponding product system over $\NN^k$ as in
\cite{RS1}. Then $A = c_0(\Lambda^0)$, and for $n \in \NN^k$, $I_n =
\clsp\{\delta_v : v \in \Lambda^0, v\Lambda^{m} = \emptyset\text{ for
all } m \le n\}$. For $m \le n \in \NN^k$, we therefore have
\[
X_m \cdot I_{n-m} = \clsp\{\delta_\mu : \mu \in \Lambda^m,
s(\mu)\Lambda^{e_i} = \emptyset\text{ whenever } m_i < n_i\}.
\]
In the language of \cite{RSY1, RSY2}, this spanning set is
familiar: $\delta_\mu \in X_m \cdot I_{n-m}$ if and only if
$\mu \in \Lambda^m \cap \Lambda^{\le n}$. That is, as a vector
space, $\Xle{n} = c_0(\Lambda^{\le n})$.
\end{example}

We thank Sean Vittadello for pointing out a simplification of
the proof of the following Lemma.

\begin{lemma} \label{lem:invariant}
Let $(G,P)$ be a quasi-lattice ordered group, let $X$ be a
compactly aligned product system over $P$ of right-Hilbert
$A$--$A$ bimodules, and fix $p,q,r \in P$ such that $pr \le q$
and $p \not= e$. Then $X_{pr} \cdot I_{(pr)^{-1}q}$ is
invariant under $\iota_{p}^{pr}(S)$ for all $S \in \Ll(X_p)$.
\end{lemma}
\begin{proof}
If $x \in X_p$, $y \in X_r$ and $i \in I_{(pr)^{-1}q}$, then
\[
\iota_{p}^{pr}(S)(xy\cdot i) = (Sx)(y\cdot i)
= ((Sx)y)\cdot i\in X_{pr} \cdot I_{(pr)^{-1}q}
\]
by definition of $\iota_{p}^{pr}$, axiom~(3) for product
systems, and associativity of multiplication in $X$. Since
vectors of the form $xy \cdot i$ where $x \in X_p$, $y \in X_r$
and $i \in I_{(pr)^{-1}q}$ span a dense subspace of $X_{pr}
\cdot I_{(pr)^{-1}q}$, and since $\iota_{p}^{pr}(S) \in
\Ll(X_{pr})$ is bounded and linear, the result follows.
\end{proof}

\begin{rmk}
Since $\iota^r_e = \phi_r$, each $X_{r} \cdot I_{r^{-1}q}$ is
also invariant under $\iota^r_e(a)$ for all $a \in \Kk(X_e)$.
\end{rmk}

\begin{notation}
Let $(G,P)$ be a quasi-lattice ordered group, and let $X$ be a
compactly aligned product system over $P$ of right-Hilbert
$A$--$A$ bimodules. Recall that $\iota^q_p$ is the zero
homomorphism from $\Ll(X_p)$ to $\Ll(X_q)$ when $p \not\le q$.
Lemma~\ref{lem:invariant} implies that for all $p, q \in P$
with $p \not= e$ there is a homomorphism from $\Ll(X_p)$ to
$\Ll(\Xle{q})$ determined by
\[
S \mapsto \bigoplus_{r \le q} \iota^r_p(S)\quad
\text{for all $S \in \Ll(X_p)$.}
\]
We denote this homomorphism by $\iotale^q_p$; it is characterised by
$(\iotale^q_p(S)x)(r) = \iota^r_p(S) x(r)$. As with the $\iota^q_p$,
when $p = e$, we write $\iotale^p_e$ for the homomorphism from
$\Kk(X_e)$ to $\Ll(\Xle{q})$ obtained from $\phile_p$ and the
isomorphism $\Kk(X_e) \cong A$.
\end{notation}

\begin{rmk}
Note that $\iotale^q_p(S)$ is the zero operator on those
summands $X_r \cdot I_{r^{-1}q}$ of $\Xle{q}$ such that $p
\not\le r$. Thus $\iotale^q_p$ can alternatively be
characterised by
\[
\iotale^q_p(S) =
    \Big(\bigoplus_{r \le q, p \not \le r} 0_{\Ll(X_r \cdot I_{r^{-1}q})}\Big)
    \oplus \Big(\bigoplus_{p \le r \le q} \iota^r_p(S)|_{X_r \cdot I_{r^{-1}q}}\Big).
\]
In particular, given $S \in \Ll(X_p)$ and $q \in P$ such that
$p \not\le q$, we have $\iotale^{q}_p(S) = 0_{\Ll(\Xle{q})}$.
\end{rmk}

To formulate our Cuntz-Pimsner covariance condition, we require
another definition.

\begin{dfn}\label{dfn:large s}
Let $(G,P)$ be a quasi-lattice ordered group. We say that a
predicate statement $\mathcal{P}(s)$ (where $s \in P$) is true
\emph{for large $s$} if: for every $q \in P$ there exists $r
\in P$ such that $q \le r$ and $\mathcal{P}(s)$ holds for all
$s \ge r$.
\end{dfn}

We now present our definition of Cuntz-Pimsner covariance. We
give a definition only in the situation that the homomorphisms
$\phile_q : A \to \Ll(\Xle{q})$ are all injective. We will see
in Lemma~\ref{lem:phi_le injective} that this is automatically
true for extensive classes of product systems. We will also see
in Example~\ref{eg:counterexample} that this is a necessary
assumption for our definition to satisfy criterion~(A) of
Section~\ref{sec:intro ps}.

\begin{dfn}\label{dfn:CP covariance}
Let $(G,P)$ be a quasi-lattice ordered group, and let $X$ be a
compactly aligned product system over $P$ of right-Hilbert
$A$--$A$ bimodules in which the homomorphisms $\phile_q$ of
Definition~\ref{dfn:Xle def} are all injective. Let $\psi$ be a
representation of $X$ in a $C^*$-algebra $B$. We say that
$\psi$ is \emph{Cuntz-Pimsner covariant} if
\begin{itemize}
\item[(${\rm CP}$)] $\sum_{p \in F} \psi^{(p)}(T_p) = 0_B$ for
    every finite $F \subset P$, and every choice of generalised
    compact operators $\{T_p \in \Kk(X_p) : p \in F\}$ such that
    $\sum_{p \in F} \iotale_p^s(T_p) = 0$ for large $s$.
\end{itemize}
\end{dfn}

\begin{rmk}
The idea is that as the indices $q \in P$ become arbitrarily
large, the associated $\Xle{q}$ approximate a notional
``boundary" of the product system (though, at this point in
time, we know of no formal way of making this idea precise).
That is~(CP) is intended to encode relations that we would
expect to hold if we could make sense of the boundary of $X$
and let all the $\Kk(X_p)$ act on it.
\end{rmk}

Our primary object of study in this paper will be the $C^*$-algebra
which is universal for representations which are both Nica covariant
and Cuntz-Pimsner covariant.

\begin{dfn}
We shall call a representation which satisfies both
(N)~and~(CP) a \emph{Cuntz-Nica-Pimsner covariant \textup{(}or
CNP-covariant\textup{)} representation}.
\end{dfn}

\begin{prop}\label{prp:NO(X) exists}
Let $(G,P)$ be a quasi-lattice ordered group, and let $X$ be a
compactly aligned product system over $P$ of right-Hilbert
$A$--$A$ bimodules such that the homomorphisms $\phile_q$ of
Definition~\ref{dfn:Xle def} are all injective. Then there
exist a $C^*$-algebra $\NO{X}$ and a CNP-covariant
representation $j_X$ of $X$ in $\NO{X}$ such that:
\begin{enumerate}
\item \label{it:spanning} $\NO{X} = \clsp\{j_X(x) j_X(y)^*
    : x,y \in X\}$; and
\item \label{it:universal} the pair $(\NO{X}, j_X)$ is
    universal in the sense that if $\psi : X \to B$ is any
    other CNP-covariant representation of $X$, then there
    is a unique homomorphism $\intfrm{\psi} : \NO{X} \to B$
    such that $\psi = \intfrm{\psi} \circ j_X$.
\end{enumerate}
Moreover, the pair $(\NO{X}, j_X)$ is unique up to canonical
isomorphism.
\end{prop}
\begin{proof}
Let $\Tc(X)$ be the universal $C^*$-algebra generated by a Nica
covariant representation $i_X$ of $X$ as in \cite{F99}. Let
$\Ii \subset \Tc(X)$ be the ideal generated by
\[\begin{split}
\Big\{\sum_{p \in F} (i_X)^{(p)}(T_p) : {}& F \subset P
\text{ is finite}, \\
& T_p \in \Kk(X_p)\text{ for each $p \in F$, and $\sum_{p \in F}
\iotale_p^s(T_p) = 0$ for large $s$}\Big\}.
\end{split}
\]
Define $\NO{X} = \Tc(X)/\Ii$, and let $q_X$ denote the quotient
map $q_X : \Tc(X) \to \NO{X}$. Let $j_X : X \to \NO{X}$ denote
the composition $j_X = q_X \circ i_X$.

Since $i_X$ is a Nica covariant representation of $X$ and $q_X$
is a homomorphism, $j_X$ satisfies $({\rm T1})$--$({\rm T3})$
and~$({\rm N})$. The definition of $q_X$ ensures that $j_X$
satisfies~(CP). Hence $j_X$ is a Cuntz-Pimsner covariant
representation of $X$. Statement~\eqref{it:spanning} follows
from the same identity for $\Tc(X)$ \cite[Equation~(6.1)]{F99}.
If $\psi : X \to B$ is an CNP-covariant representation then it
is, in particular, a Nica covariant representation of $X$ and
it follows from \cite[Theorem~6.3]{F99} that there is a
homomorphism $\psi_* : \Tc(X) \to B$ such that $\psi = \psi_*
\circ i_X$. Since $\psi$ is Cuntz-Pimsner covariant, we have
$\Ii \subset \ker(\psi_*)$, and it follows that $\psi_*$
descends to a homomorphism $\intfrm{\psi} : \NO{X} \to B$,
establishing~\eqref{it:universal}.

All that remains to be proved is the uniqueness claim, for
which we give the following standard argument. If $(C, \rho)$
is another such pair then~(2) for $\NO{X}$ implies that there
is a homomorphism $\intfrm{\rho} : \NO{X} \to C$ such that
$\intfrm{\rho}(j_X(x)) = \rho(x)$ for all $x$. Statement~(2)
for $C$ implies that there is a homomorphism $\intfrm{j_X} : C
\to \NO{X}$ such that $\intfrm{j_X}(\rho(x)) = j_X(x)$ for all
$x$. Applications of~(1) then show that $\intfrm{\rho}$ and
$\intfrm{j_X}$ are surjective and are mutually inverse.
\end{proof}

\begin{rmk}
No obvious notation for the homomorphism $\intfrm{\psi}$ of
Proposition~\ref{prp:NO(X) exists}(\ref{it:universal}) occurred
to us. We were loathe to re-define Fowler's notation $\psi_*$:
our $\intfrm{\psi}$ is induced from Fowler's $\psi_*$ by
regarding $\NO{X}$ as a quotient of $\Tc(X)$, so employing the
same notation would lead to confusion in any situation where
representations of both $\Tc(X)$ and $\NO{X}$ are discussed. We
settled on $\intfrm{\psi}$ on the basis that it might bring to
mind the integrated form $\pi \times \psi$ of a representation
$(\pi,\psi)$ of a single bimodule.
\end{rmk}

\begin{rmk}
It should be emphasised that $\NO{X}$ may differ from Fowler's
$\Oo_X$ even should our notion of Cuntz-Pimsner covariance and
Fowler's coincide. The algebras $\NO{X}$ and $\Oo_X$ will only
coincide for product systems such that: the two versions of
Cuntz-Pimsner covariance are equivalent; and Cuntz-Pimsner
covariance implies Nica covariance.
\end{rmk}

We conclude this section by showing that the hypothesis that
the $\phile_q$ are all injective is automatic for broad classes
of product systems; but we also show by example that there
exist product systems for which this hypothesis fails.

Specifically, we shall show that the $\phile_q$ are all
injective whenever the $\phi_q$ are all injective, and also
whenever the quasi-lattice ordered pair $(G,P)$ has the
property that every nonempty bounded subset of $P$ contains a
maximal element in the following sense.
\begin{equation}\label{eq:max elements}
\parbox{0.8\textwidth}{
If $S \subset P$ is nonempty and there exists $q \in P$ such that
$p \le q$ for all $p \in S$, then there exists $p \in S$ such that
$p \not\le p'$ for all $p' \in S \setminus \{p\}$.}
\end{equation}

\begin{lemma}\label{lem:phi_le injective}
Let $(G,P)$ be a quasi-lattice ordered group, and let $X$ be a
compactly aligned product system over $P$ of right-Hilbert
$A$--$A$ bimodules. If each $\phi_p$ is injective, or if $P$
satisfies~\eqref{eq:max elements}, then $\phile_{q} : A \to
\Ll(\Xle{q})$ is injective for each $q \in P$.
\end{lemma}
\begin{proof}
Suppose first that each $\phi_p$ is injective. Then $\Xle{p}
\cong X_p$ for all $p$ and this isomorphism intertwines
$\phile_p$ and $\phi_p$, so the $\phile_p$ are injective.

Now suppose that $P$ satisfies~\eqref{eq:max elements}. Fix $a
\in A \setminus \{0\}$ and $q \in P$. We must show that
$\phile_{q}(a) \not= 0$; that is, $\phi_p(a)|_{X_p \cdot
I_{p^{-1} q}} \not= 0$ for some $p \le q$. We have $\phi_e(a)
\not= 0$ because $\phi_e(a)(a^*) = a a^* \not= 0$. Let $S = \{p
\in P : p \le q, \phi_p(a) \not= 0\}$. Then $S$ is nonempty,
and is bounded above by $q$. Since $P$ satisfies~\eqref{eq:max
elements}, it follows that $S$ contains a maximal element $p$.
Since $p \in S$, we have $\phi_p(a) \not=0$, so we may fix $x
\in X_p$ such that $\phi_p(a) x \not= 0$.

Since $p$ is maximal in $S$, if $e < r \le p^{-1}q$ then $a \in
\ker(\phi_{pr})$. So if $y \in X_r$ where $e < r \le p^{-1}q$,
then $(\phi_p(a) x)y = \phi_{pr}(a)(xy) = 0$; that is,
$(\phi_p(a) x) y = 0$ for all $e < r \le p^{-1}q$ and all $y
\in X_r$. Lemma~\ref{lem:Xle characterisation} therefore
implies that $\phi_p(a)x \in X_p \cdot I_{p^{-1}q}$.

Let $\{\mu_\lambda : \lambda \in \Lambda\}$ be an approximate
identity for $I_{p^{-1} q}$. By the preceding paragraph,
$\{\phi_p(a) x \cdot \mu_\lambda\}_{\lambda \in \Lambda}$ is
norm-convergent to $\phi_p(a) x \not= 0$, so there exists
$\lambda \in \Lambda$ such that $\phi_p(a) x \cdot \mu_\lambda
\not= 0$. Setting $y = x \cdot \mu_\lambda$, we have $y \in X_p
\cdot I_{p^{-1} q}$ and $\phi_p(a) y \not=0$. That is
$\phi_p(a)|_{X_p \cdot I_{p^{-1}q}} \not= 0$, so $\phile_{q}(a)
\not= 0$ by definition.
\end{proof}

The hypotheses of Lemma~\ref{lem:phi_le injective} are not just
an artifact of our proof. To see why, consider the following
example.

\begin{example}\label{eg:counterexample}
Let $G = \ZZ \times \ZZ$, and let $P = ((\NN\setminus\{0\})
\times \ZZ) \cup (\{0\} \times \NN)$. Then $P$ is a
subsemigroup of $G$ satisfying $P \cap -P = \{(0,0)\}$. The
partial order $\le$ on $G$ defined by $m \le n \iff n-m \in P$
is the lexicographic order on $\ZZ \times \ZZ$, and in
particular $(G,P)$ is a quasi-lattice ordered group. Note that
$(G,P)$ does not satisfy~\eqref{eq:max elements}: let $S$
denote the subset $S = \{0\} \times \NN$. Then $S$ is bounded
above by $(1,0)$, but has no maximal element.

Consider the right-Hilbert $\CC^2$ module $\CC^2_{\CC^2}$ with
the usual right action and inner-product. Then
$\Ll(\CC^2_{\CC^2}) = \Kk(\CC^2_{\CC^2}) \cong \CC^2$: the
element $(z_1, z_2) \in \CC^2$ acts by point-wise left
multiplication on $\CC^2_{\CC^2}$. Define homomorphisms
$\phi_{\Sbit}, \phi_{\notSbit} : \CC^2 \to \Ll(\CC^2_{\CC^2})$
by $\phi_{\Sbit} = \id_{\CC^2}$ and $\phi_{\notSbit}(z_1, z_2)
= (z_1, z_1)$. We write $X_{\Sbit}$ (respectively
$X_{\notSbit}$) for $\CC^2_{\CC^2}$ regarded as a right-Hilbert
$\CC^2$--$\CC^2$ bimodule with left action implemented by
$\phi_{\Sbit}$ (respectively $\phi_{\notSbit}$).

There are isomorphisms
\begin{align*}
X_{\Sbit} \otimes_{\CC^2} X_{\Sbit} &\cong X_{\Sbit}
&&\text{ determined by } (z_1, z_2) \otimes_{\CC^2} (w_1, w_2) \mapsto (z_1w_1, z_2w_2), \\
X_{\Sbit} \otimes_{\CC^2} X_{\notSbit} &\cong X_{\notSbit}
&&\text{ determined by } (z_1, z_2) \otimes_{\CC^2} (w_1, w_2) \mapsto (z_1w_1, z_1w_2), \\
X_{\notSbit} \otimes_{\CC^2} X_{\Sbit} &\cong X_{\notSbit}
&&\text{ determined by } (z_1, z_2) \otimes_{\CC^2} (w_1, w_2) \mapsto (z_1w_1, z_2w_2),\text{ and} \\
X_{\notSbit} \otimes_{\CC^2} X_{\notSbit} &\cong X_{\notSbit}
&&\text{ determined by } (z_1, z_2) \otimes_{\CC^2} (w_1, w_2) \mapsto (z_1w_1, z_1w_2);
\end{align*}
to see this, one checks that each of these formulae preserves
inner-products of elementary tensors.

For $p \in S$, let $X_p = X_{\Sbit}$, and for $p \in P
\setminus S$, let $X_p = X_{\notSbit}$. With multiplication
maps defined as above, $X$ is a product system over $P$. Note
that for $p \in P$, $\phi_p = \phi_{\Sbit}$ if $p \in S$, and
$\phi_p = \phi_{\notSbit}$ if $p \not\in S$. Since $\Ll(X_p) =
\Kk(X_p)$ for all $p$, the left action of $\CC^2$ on each fibre
is by compact operators, and in particular $X$ is compactly
aligned.

We have $(0,1) \le p$ for all $p \in P \setminus\{(0,0)\}$.
Since $\ker(\phi_{(0,1)}) = \{0\}$, it follows that $I_p =
\{0\}$ for $p \not= (0,0)$. Hence $\Xle{q} = X_q$ and $\phile_q
= \phi_q$. In particular, $\phile_{(1,0)} = \phi_{(1,0)} =
\phi_{\notSbit}$ is not injective.

Observe that $\iotale^{s}_{(0,0)}((0,1)) = 0_{\Ll(\Xle{s})}$
for large $s$. Hence every representation $\psi$ of $X$
satisfying~(CP) satisfies $\psi_e((0,1)) = 0$. In particular
the algebra universal for such representations does not satisfy
criterion~(A) of Section~\ref{sec:intro ps}.
\end{example}

\begin{rmk}
Subject to failure of the hypothesis of Lemma~\ref{lem:phi_le
injective}, Example~\ref{eg:counterexample} is as well-behaved
as possible: $(G,P)$ is countable and totally ordered, and the
natural order topology is discrete; $A$ and the $X_p$ are all
finite-dimensional, so that in particular the action on each
fibre is by compact operators; and each $X_p$ is essential as a
left $A$-module in the sense that $\phi_p(A)X_p = X_p$.
\end{rmk}

\section{Injectivity of the universal CNP-covariant representation}

In this section we prove our main theorem.

\begin{theorem}\label{thm:isometric}
Let $(G,P)$ be a quasi-lattice ordered group, and let $X$ be a
compactly aligned product system over $P$ of right-Hilbert
$A$--$A$ bimodules. Suppose that the homomorphisms $\phile_q$
of Definition~\ref{dfn:Xle def} are all injective. Then the
universal CNP-covariant representation $j_X : X \to \NO{X}$ is
isometric: $\|j_X(x)\| = \|x\|$ for all $x \in X$. In
particular, the conclusion holde if each $\phi_p$ is injective,
or if $P$ satisfies~\eqref{eq:max elements}.
\end{theorem}

We now introduce a modification, based on the $\Xle{p}$, of
Fowler's Fock representation.

\begin{notation}
Let $(G,P)$ be a quasi-lattice ordered group, and let $X$ be a
product system over $P$ of right-Hilbert $A$--$A$ bimodules. As
on \cite[page~340]{F99}, we let $F(X) = \bigoplus_{p \in P}
X_p$, and call it the \emph{Fock space} of $X$. We also define
$\aF(X) = \bigoplus_{p \in P} \Xle{p}$, and call it the
\emph{augmented Fock space} of $X$.
\end{notation}

Fowler shows \cite[page~340]{F99} that for $x \in X_p$ there is an
adjointable operator $l(x)$ on $F(X)$ determined by
\[
(l(x)y)(q) =
\begin{cases}
    x (y(p^{-1}q)) &\text{ if $p \le q$}\\
    0 &\text{ otherwise.}
\end{cases}
\]
He shows further that $l : X \to \Ll(F(X))$ is a Nica covariant
representation of $X$. The next lemma shows that we obtain a parallel
result if we replace $F(X)$ with $\aF(X)$.

\begin{lemma}\label{lem:Augl}
Let $(G,P)$ be a quasi-lattice ordered group, and let $X$ be a
compactly aligned product system over $P$ of right-Hilbert
$A$--$A$ bimodules.
\begin{enumerate}
\item Let $p,q,r \in P$, with $q \le r$, and let $x \in
    X_p$ and $z \in X_q \cdot I_{q^{-1} r}$. Then $xz \in
    X_{pq} \cdot I_{(pq)^{-1}pr}$.
\item\label{it:def ltilde} For $p \in P$ and $x \in X_p$,
    there is an adjointable operator $\Augl(x) \in
    \Ll(\aF(X))$ which satisfies
\[
(\Augl(x)y)(q) =
\begin{cases}
x(y(p^{-1}q)) &\text{ if $p \le q$} \\
0 &\text{ otherwise}
\end{cases}
\]
for all $y \in \aF(X)$. Moreover if $y \in \aF(X)$
satisfies $y(s) = 0$ for all $s \ge p$, then $\Augl(x)^*y =
0$.
\item The map $x \mapsto \Augl(x)$ of~$($\ref{it:def
    ltilde}$)$ is a Nica covariant representation of $X$ in
    $\Ll(\aF(X))$.
\end{enumerate}
\end{lemma}
\begin{proof}
To prove~(1), write $z = z' \cdot i$ where $z' \in X_q$ and $i
\in I_{q^{-1}r} \subset A = X_e$. By associativity of
multiplication in $X$ and condition~(4) of the definition of a
product system we have $xz = x(z' \cdot i) = (xz')\cdot i \in
X_{pq} \cdot I_{(pq)^{-1}(pr)}$.

For~(2), observe that each $\Xle{q}$ is a sub-module of $F(X)$,
and that the restriction of $l(x)$ to this submodule agrees
with $\Augl(x)$. Since $l(x)$ is an adjointable operator on
$F(X)$, the first part of~(2) follows. For the second part,
observe that by linearity and continuity of $\Augl(x)^*$, we
may assume that $y$ has just one nonzero coordinate. That is,
we may assume that there exists $r \in P$ such that $r \not\ge
p$ and $y(s) = 0$ for $s \not= r$. For $z \in \aF(X)$, we then
have
\[
\langle \Augl(x)^* y, z \rangle_A
 = \langle y, \Augl(x) z \rangle_A
 = \sum_{s \in P} \langle y(s), (\Augl(x) z)(s) \rangle^{\Xle{s}}_A
 = \langle y(r), (\Augl(x) z)(r) \rangle^{\Xle{r}}_A
 = 0
\]
by the first part of~(2).

For~(3), recall that $l$ satisfies (T1)--(T3)~and~(N). It
follows that $\Augl$ does as well.
\end{proof}

Let $B$ be the $C^*$-algebra generated by $\{\Augl(x) : x \in X\}
\subset \Ll(\aF(X))$. Let $\Ii_{\aF(X)}$ be the ideal in $B$
generated by
\begin{align*}
\Big\{ \sum_{p \in F} \Augl^{(p)}(T_p) : {}&
F \subset P \text{ is finite}, \\
&T_p \in \Kk(X_p)\text{ for each $p \in F$, and $\sum_{p \in F}
\iotale_p^s(T_p) = 0$ for large $s$}\Big\}.
\end{align*}
Let $q_{\aF(X)} : B \to B/\Ii_{\aF(X)}$ denote the quotient
map. Then $\psi := q_{\aF(X)} \circ \Augl$ is a CNP-covariant
representation of $X$ in $B/\Ii_{F(X)}$.

\begin{prop} \label{prp:faithful on A}
Let $(G,P)$ be a quasi-lattice ordered group, and let $X$ be a
compactly aligned product system over $P$ of right-Hilbert
$A$--$A$ bimodules. Suppose that each $\phile_q$ is injective.
Then $q_{F(X)} \circ \Augl : X \to B/\Ii_{\Ff(X)}$ is faithful
on $A$.
\end{prop}
\begin{proof}
We must show that $\Ii_{\aF(X)} \cap \Augl(A) = \{0\}$.

For $a \in A$, the restriction $\Augl(a)|_{\Xle{q}}$ of
$\Augl(a)$ to the $q$-summand of $\aF(X)$ is given by
$\Augl(a)x = \phile_{q}(a) x$. Since $\phile_q$ is injective
and hence isometric, it follows that the operator norm
$\|\Augl(a)|_{\Xle{q}}\|$ is equal to the $C^*$-norm $\|a\|$ of
$a \in A$. It therefore suffices to show that
\begin{equation}\label{eq:vanishing at boundary}
\text{for each $\psi \in \Ii_{\aF(X)}$ and each $\varepsilon > 0$
there exists $s \in P$ such that $\|\psi|_{\Xle{s}}\| <
\varepsilon$.}
\end{equation}

We do this in stages. Let $K \subset \Ii_{\aF(X)}$ be the subset
\[\begin{split}
K = \Big\{ \sum_{p \in F} \Augl^{(p)}(T_p) : {}
& F \subset P\text{ is finite}, \\
&T_p \in \Kk(X_p)\text{ for each $p \in F$, and $\sum_{p \in F}
\iotale_p^s(T_p) = 0$ for large $s$}\Big\}.
\end{split}
\]
Then $K$ is a subspace of $\Ii_{\aF(X)}$, and a continuity argument
shows that for each $\psi \in \overline{K}$, each $\varepsilon
> 0$, and each $q \in P$ there exists $r \ge q$ such that
$\|\psi|_{\Xle{s}}\| < \varepsilon$ for $s \ge r$.

By \cite[Proposition~5.10]{F99}, we have $B = \clsp\{\Augl(x)
\Augl(y)^* : x,y \in X\}$. Consequently,
\[
\Ii_{\aF(X)} = \clsp\{\Augl(x)\Augl(y)^* k \Augl(x') \Augl(y')^* : k\in
\overline{K}, x,y,x',y', \in X\}.
\]
Fix $k \in \overline{K}$ and $x,y,x',y', \in X$ --- say $x \in
X_{p(x)}$, $y \in X_{p(y)}$, $x' \in X_{p(x')}$ and $y' \in
X_{p(y')}$ ---  and fix $q \in P$. We claim that there exists
$r \ge q$ such that $\|\Augl(x)\Augl(y)^* k \Augl(x')
\Augl(y')^*|_{\Xle{s}}\| < \varepsilon$ for each $s \ge r$. To
see this, suppose first that $p(y') \vee q < \infty$. Then
there exists $r' \ge p(x')p(y')^{-1}(p(y') \vee q)$ such that
$\|k|_{\Xle{s'}}\| <
\frac{\varepsilon}{\|x\|\|y\|\|x'\|\|y'\|}$ for all $s' \ge
r'$. Let $r = p(y')p(x')^{-1}r' \ge q$. Then for each $s \ge
r$, we have $p(x')p(y')^{-1} s \ge r'$, and hence
\[
\|\Augl(x)\Augl(y)^* k \Augl(x') \Augl(y')^*|_{\Xle{s}}\| < \varepsilon.
\]
Now suppose that $p(y') \vee q = \infty$. Then $s \ge q$
implies $s \not\ge p(y')$, so Lemma~\ref{lem:Augl}(2) implies
that $\Augl(y')^*|_{\Xle{s}} = 0$. Consequently $r = q$
satisfies $\Augl(x)\Augl(y)^* k \Augl(x')
\Augl(y')^*|_{\Xle{s}} = 0$ for every $s \ge r$.

The assertion~\eqref{eq:vanishing at boundary} now follows by
linearity and continuity.
\end{proof}

\begin{rmk}
Since each summand of $\aF(X)$ is a sub-module of $F(X)$, it
was not necessary to introduce $\aF(X)$ at all. Essentially the
argument of Proposition~\ref{prp:faithful on A} is valid if we
work instead with Fowler's $l : X \to F(X)$ followed by the
appropriate quotient map. However using $\aF(X)$ makes the
argument clearer, and helps give some intuition for the
significance of the $\Xle{p}$.
\end{rmk}

\begin{proof}[Proof of Theorem~\ref{thm:isometric}]
Proposition~\ref{prp:faithful on A} together with he universal
property of $\NO{X}$ implies that $j_X$ is injective on $X_{e}
= A$. As $j_X|_{X_e}$ is a $C^*$-homomorphism, it follows that
$j_X$ is isometric on $A$. For any $x \in X$, we therefore have
\[
\|x\|^2 = \|\langle x,x \rangle_A \| = \|j_X(\langle x,x\rangle_A)\| =
\|j_X(x)^* j_X(x)\| = \|j_X(x)\|^2,
\]
completing the proof.
\end{proof}

\begin{rmk}
Note that $\NN^k$, and indeed every right-angled Artin
semigroup, satisfies~\eqref{eq:max elements}, so
Theorem~\ref{thm:isometric} applies when $P$ is any of these
semigroups.
\end{rmk}

\begin{rmk}
Theorem~\ref{thm:isometric} shows that if each $\phile_q$ is
injective, and in particular if $P$ satisfies~\eqref{eq:max
elements}, then $\NO{X}$ satisfies criterion~(A) of
Section~\ref{sec:intro ps}. In the next section, we will show
that it also satisfies~(B) in a number of motivating examples.
We will achieve this indirectly by combining
Theorem~\ref{thm:isometric} with the uniqueness theorems
of~\cite{CL, K6, RSY2}.

We have not verified~(B). This is done in~\cite{CLSV} for
certain pairs $(G,P)$ using a careful analysis of the
fixed-point algebra $\fixNO{X}{\delta} = \clsp\{j_X^{(p)}(T_p)
: p \in P, T_p \in \Kk(X_p)\}$.
\end{rmk}

\section{Relationships to other constructions}\label{sec:comparison}

In this section we discuss the relationship between $\NO{X}$
and a number of other $C^*$-algebras. We begin by showing that
when each pair in $P$ has a least upper bound and each $\phi_p$
is injective and takes values in $\Kk(X_p)$, our $\NO{X}$
coincides with Fowler's Cuntz-Pimsner algebra
$\Oo_X$~\cite{F99}. We then demonstrate that our $\NO{X}$ is
compatible with Katsura's $\Oo_X$ \cite{K5,K6}, with the
Cuntz-Krieger algebras of finitely aligned higher-rank
graphs~\cite{RSY2}, and with Crisp and Laca's boundary
quotients of Toeplitz algebras~\cite{CL}.

\subsection{Fowler's Cuntz-Pimsner algebras}\label{sec:Fowler's}

Kumjian and Pask's uniqueness theorems for higher-rank graph
$C^*$-algebras \cite{KP} suggest that Fowler's notion of
Cuntz-Pimsner covariance determines a universal $C^*$-algebra
which satisfies (A)~and~(B) of Section~\ref{sec:intro ps} when
$(G,P) = (\ZZ^k, \NN^k)$, and each $\phi_p$ is injective with
$\phi_p(A) \subset \Kk(X_p)$. Proposition~5.4 of~\cite{F99}
then suggests that we may be able to relax the requirement that
$(G,P) = (\ZZ^k, \NN^k)$ and insist only that each pair in $P$
has a common upper bound. The next proposition shows that in
this case~$\NO{X}$ and Fowler's~$\Oo_X$ coincide.

Note that we assume only that each pair in $P$ has a common
upper bound, not that each pair in $G$ has a common upper
bound. Under the latter hypothesis, we deduce that each pair in
$G$ has a least common upper bound, and hence that $G$ is in
fact lattice ordered (take $g \wedge h = (g^{-1} \vee
h^{-1})^{-1})$. It is not clear to us whether the assumption
that each pair in $P$ has a common upper bound implies that $P$
is lattice-ordered, or even that $G$ is.

Recall that a representation $\psi$ of $X$ is Cuntz-Pimsner
covariant in the sense of \cite[Definition~2.5]{F99} if,
whenever $\phi_p(a) \in \Kk(X_p)$, we have $\psi^{(p)}(\phi(a))
= \psi_e(a)$.

\begin{prop}\label{prop:ours->Fowlers}
Let $(G,P)$ be a quasi-lattice ordered group, and let $X$ be a
compactly aligned product system over $P$ of right-Hilbert $A$--$A$
bimodules. Suppose that each pair in $P$ has a least upper bound.
Suppose that for each $p \in P$, the homomorphism $\phi_p : A \to
\Ll(X_p)$ is injective. Let $\psi : X \to B$ be a representation of
$X$.
\begin{enumerate}
\item\label{it:if} If $\psi$
    satisfies~\textup{(}CP\textup{)}, then it is
    Cuntz-Pimsner covariant in the sense of
    \cite[Definition~2.5]{F99}.
\item\label{it:onlyif} If $\phi_p(A) \subset \Kk(X_p)$ for
    each $p \in P$, and $\psi$ is Cuntz-Pimsner covariant
    in the sense of \cite[Definition~2.5]{F99}, then $\psi$
    satisfies~\textup{(}CP\textup{)}.
\end{enumerate}
\end{prop}
\begin{proof}
We begin with some observations which we will use to prove both
\eqref{it:if}~and~\eqref{it:onlyif}. Let $a \in A$. Recall that
$X_e$ is equal to $A$, and that $\phi_e(a) : X_e \to X_e$ is
the left-multiplication operator $b \mapsto ab$, which we
denote by $L_a \in \Kk(X_e)$. Every element $T$ of $\Kk(X_e)$
can be written as $T = L_a$ for some $a \in A$ (see
\cite[Lemma~2.26]{TFB}). By definition of $\iota^p_e$, we have
\begin{equation}\label{eq:lifting a}
\iota^p_e(L_a) = \phi_p(a)\quad\text{ for all $p \in P$, $a \in A$.}
\end{equation}

We now prove \eqref{it:if}~and~\eqref{it:onlyif} separately.

For~(\ref{it:if}), let $p \in P$ and $a \in A$, and suppose
that $\phi_p(a) \in \Kk(X_p)$. We must show that $\psi_e(a) -
\psi^{(p)}(\phi_p(a)) = 0$. Equation~\eqref{eq:lifting a}
implies that $\iota^{s}_e(L_a) - \iota^{s}_p(\phi_p(a)) =
0_{\Ll(X_{s})}$ for all $s \ge p$. Since each pair in $P$ has a
least upper bound, it follows that for each $q \in P$, we have
$\iota^{s}_e(L_a) - \iota^{s}_p(\phi_p(a)) = 0_{\Ll(X_{s})}$
for all $s \ge p \vee q$.

Since each $\phi_p$ is injective, we have $I_p = \{0\}$ for $p
\not=e$. Hence $\Xle{p} = X_p$ and $\iotale^q_p = \iota^q_p$ for all
$p \le q \in P$. Thus $\iotale_e^s(L_a) - \iotale_p^s(\phi_p(a)) = 0$
for large $s$. Since $\psi$ satisfies~(CP), we deduce that
$\psi^{(e)}(L_a) - \psi^{(p)}(\phi_p(a)) = 0$. A straightforward
computation using an approximate identity for $A$ and that $\psi_e$
is a homomorphism shows that $\psi^{(e)}(L_a) = \psi_e(a)$. Hence
$\psi_e(a) = \psi^{(p)}(\phi_p(a))$, and $\psi$ is Cuntz-Pimsner
covariant in the sense of \cite[Definition~2.5]{F99}.

Now for~(\ref{it:onlyif}), fix a finite subset $F \subset P$
and compact operators $T_p \in \Kk(X_p)$, $p \in F$ such that
$\sum_{p \in F} \iotale_p^s(T_p) = 0$ for large $s$. An
inductive argument using that every pair in $P$ has a least
upper bound shows that there is a least upper bound $q =
\bigvee F$ for $F$ in $P$. Since each $\phi_p$ is injective, we
have $\iota^s_q$ injective for $s \ge q$, so we must have
$\sum_{p \in F} \iota^q_p(T_p) = 0$.

We may assume without loss of generality that $e \in F$, and
rearrange to obtain $\sum_{p \in F \setminus \{e\}}
\iota^q_p(T_p) = -\iota^q_e(T_e)$. As observed above, $T_e =
L_a$ for some $a \in A$, so~\eqref{eq:lifting a} implies that
$\iota^q_e(T_e) = \phi_q(a)$. Since the left action of $A$ on
$X_q$ is by compact operators, we have $\phi_q(a) \in
\Kk(X_q)$, and hence Fowler's Cuntz-Pimsner covariance
condition forces
\begin{equation}\label{eq:take a across}
\psi^{(q)}\Big(\sum_{p \in F \setminus \{e\}} \iota^q_p(T_p)\Big) =
-\psi_e(a).
\end{equation}
Since $A$ acts compactly on the left of each $X_p$,
\cite[Corollary~3.7]{Pim97} shows that each $\iota^q_p(T_p) \in
\Kk(X_q)$, and the argument of \cite[Lemma~3.10]{Pim97} shows that
$\psi^{(q)}(\iota^q_p(S)) = \psi^{(p)}(S)$. Hence~\eqref{eq:take a
across} implies that
\[
\sum_{p \in F \setminus \{e\}} \psi^{(p)}(T_p) = -\psi_e(a).
\]
Since $\psi_e(a) = \psi^{(e)}(L_a) = \psi^{(e)}(T_e)$, we have
$\sum_{p \in F} \psi^{(p)}(T_p) = 0$. So $\psi$ satisfies~(CP).
\end{proof}

\begin{cor}
Let $(G,P)$ be a quasi-lattice ordered group, and let $X$ be a
compactly aligned product system over $P$ of right-Hilbert
$A$--$A$ bimodules. Suppose that each pair in $P$ has a least
upper bound. Suppose that each $\phi_p : A \to \Ll(X_p)$ is
injective with $\phi_p(A) \subset \Kk(X_p)$. Let $\psi : X \to
B$ be a representation of $X$. Then $\psi$ is CNP-covariant if
and only if $\psi^{(p)} \circ \phi_p = \psi_e$ for all $p \in
P$.
\end{cor}
\begin{proof}
Since each $\phi_p(A) \subset \Kk(X_p)$, the condition
$\psi^{(p)} \circ \phi_p = \psi_e$ for all $p \in P$ is
precisely the Cuntz-Pimsner covariance condition of
\cite[Definition~2.5]{F99}. Proposition~5.4 of \cite{F99}
implies that if $\psi$ is Cuntz-Pimsner covariant in this
sense, then it is also Nica covariant. Hence the result follows
from Proposition~\ref{prop:ours->Fowlers}.
\end{proof}

\subsection{Katsura's $C^*$-algebras associated to Hilbert bimodules}\label{sec:Katsura's}

Let $X$ be a right-Hilbert $A$--$A$ bimodule. There is a
product system $X^\otimes$ over $\NN$ of right-Hilbert $A$--$A$
bimodules such that $X^\otimes_0 = A$, and $X^\otimes_n =
X^{\otimes n}$ for $1 \le n \in \NN$. For nonzero $m,n$, the
isomorphism $M_{m,n} : X^\otimes_m \otimes X^\otimes_n \to
X^\otimes_{m+n}$ implementing the multiplication in the system
is the natural isomorphism $X^{\otimes m} \otimes X^{\otimes n}
\cong X^{\otimes m+n}$. As $\NN$ is totally ordered,
$X^\otimes$ is compactly aligned.

If $(\psi,\pi)$ is a representation of $X$, then there is a
representation $\psi^\otimes$ of $X^\otimes$ given by
$\psi^\otimes_e = \pi$ and $\psi^\otimes_n = \psi^{\otimes n}$
for $n \ge 1$. This representation is automatically Nica
covariant, and every representation of $X^\otimes$ is of this
form. Proposition~2.11 of \cite{F99} says that if the
homomorphism $\phi : A \to X$ implementing the left action is
injective or takes values in $\Kk(X)$ then $\psi$ is covariant
in Pimsner's sense if and only if $\psi^{\otimes}$ is
Cuntz-Pimsner covariant in Fowler's sense. Hence Pimsner's
$\Oo_X$ \cite{Pim97} and Fowler's $\Oo_{X^\otimes}$ \cite{F99}
coincide.

A key goal of our construction was to achieve the same outcome
with respect to Katsura's reformulation of Cuntz-Pimsner
covariance for a single Hilbert bimodule. That is, given an
arbitrary Hilbert bimodule $X$, we desire that
$\NO{X^{\otimes}}$ should coincide with Katsura's
$\Oo_X$~\cite{K5, K6}.

Recall that if $I \lhd A$ is an ideal in a $C^*$-algebra, then
$I^\perp$ denotes the ideal $\{a \in A : ab = 0\text{ for all }
b \in I\}$. Recall also that Katsura's $\Oo_X$ is the universal
$C^*$-algebra generated by a representation $(i_A,i_X)$ of $X$
such that $i_X^{(1)}(\phi(a)) = i_A(a)$ for all $a \in
\ker(\phi)^\perp$ such that $\phi(a) \in \Kk(X)$.

\begin{prop}\label{prp:same as Katsie's}
Let $X$ be a right-Hilbert $A$--$A$ bimodule. Let $(i_A, i_X)$
be the universal representation of $X$ on $\Oo_X$, and
$j_{X^\otimes}$ be the universal representation of $X^\otimes$
on $\NO{X^\otimes}$.
\begin{itemize}
\item[(1)] There is an isomorphism $\theta : \Oo_X \to
    \NO{X^\otimes}$ satisfying $\theta(i_A(a)) =
    j_{X^\otimes}(a)$ and $\theta(i_X(x)) =
    j_{X^\otimes}(x)$ for all $a \in A$ and $x \in X$.
\item[(2)] Let $(\pi, \psi)$ be a representation of $X$ and
    let $\psi^\otimes$ be the corresponding Nica covariant
    representation of $X^\otimes$. Then $(\pi,\psi)$ is
    covariant in the sense of Katsura if and only if
    $\psi^\otimes$ satisfies~(CP).
\end{itemize}
\end{prop}
\begin{proof}
Statement~(2) follows from~(1) and the universal properties of
$\Oo_X$ and $\NO{X^\otimes}$, so it suffices to prove~(1).

We have $\ker(\phi) \subset \ker(\phi \otimes 1_{n-1})$ for $n
\ge 1$, so $I_n$ is equal to $A$ if $n = 0$ and is equal to
$\ker(\phi)$ if $n \not= 0$. Let $(j_0, j_1)$ denote the
representation of $X$ determined by the universal
representation of $X^{\otimes}$ on $\NO{X^\otimes}$. If $a \in
\ker(\phi)^\perp \cap \phi^{-1}(\Kk(X))$, then $\phile_{1}(a) =
0_{\ker(\phi)} \oplus \phi(a)$. Let $S = \phi(a) \in \Kk(X)$.
Then $\iotale^{1}_0(\phi_0(a)) - \iotale^{1}_1(S) = 0$, and it
follows that $\iotale^{n}_0(\phi_0(a)) - \iotale^{n}_1(S) = 0$
for all $n \ge 1$. Since $j_{X^\otimes}$ satisfies~(CP), we
therefore have $j_0(a) - j_1^{(1)}(S) = 0$; that is $j_0(a) =
j_1^{(1)}(\phi(a))$. Thus $(j_0, j_1)$ is covariant in the
sense of Katsura, and the universal property of $\Oo_X$ implies
that there is a homomorphism $\theta : \Oo_X \to \NO{X}$
determined by $\theta \circ i_A = j_0$ and $\theta \circ i_X =
j_1$. Moreover, $\theta$ is surjective because
\[
\NO{X}
= \clsp\{j_1(x_1) \cdots j_1(x_m) j_1(y_n)^* \cdots j_1(y_1)^*
: m,n \in \NN, x_i, y_i \in X\}.
\]
Theorem~\ref{thm:isometric} implies that $\theta \circ i_A =
j_0$ is injective, so an application of Katsura's
gauge-invariant uniqueness theorem \cite[Theorem~6.2]{K6} shows
that $\theta$ is injective.
\end{proof}

\subsection{Cuntz-Krieger algebras of finitely aligned higher-rank graphs}\label{sec:hrg}

In this section, we use the notation and conventions of
\cite{RSY2} for higher-rank graphs. In \cite{RS1}, a product
system $X(\Lambda)$ over $\NN^k$ of right-Hilbert
$c_0(\Lambda^0)$--$c_0(\Lambda^0)$ bimodules is associated to
each $k$-graph $\Lambda$. When $\Lambda$ is row-finite and has
no sources, the homomorphism $\phi_p$ implementing the left
action of $c_0(E^0)$ on $X(\Lambda)_p$ is an injective
homomorphism into the compact operators on $X(\Lambda)_p$.
Corollary~4.4 of \cite{RS1} shows that $C^*(\Lambda)$ coincides
with $\Oo_{X(\Lambda)}$ as defined in \cite{F99} and hence with
$\NO{X(\Lambda)}$ by Proposition~\ref{prop:ours->Fowlers}. A
key goal of our construction is to extend this to arbitrary
finitely aligned $k$-graphs and the corresponding compactly
aligned product systems of bimodules. We show in this section
that we have achieved this aim.

We briefly recall some salient point about $k$-graphs from
\cite{RSY2}. A $k$-graph is a countable category $\Lambda$
together with a functor $d : \Lambda \to \NN^k$ satisfying the
factorisation property: for all $m,n \in \NN^k$ and $\lambda
\in d^{-1}(m+n)$ there exist unique $\mu \in d^{-1}(m)$ and
$\nu \in d^{-1}(n)$ such that $\lambda = \mu\nu$. Each
$d^{-1}(n)$ is denoted $\Lambda^n$. The elements of $\Lambda^0$
are called vertices, and are in bijection with the objects of
$\Lambda$, so the codomain and domain maps in the category
$\Lambda$ determine maps $r,s : \Lambda \to \Lambda^0$. For
$\mu,\nu \in \Lambda$, we write $\MCE(\mu,\nu)$ for the
collection $\{\lambda \in \Lambda : d(\lambda) = d(\mu) \vee
d(\nu), \lambda = \mu\mu'=\nu\nu'\text{ for some }\mu', \nu'
\in \Lambda\}$. The $k$-graph $\Lambda$ is said to be
\emph{finitely aligned} if $\MCE(\mu,\nu)$ is finite (possibly
empty) for every pair $\mu,\nu$ of paths in $\Lambda$.

For $n \in \NN^k$, we write $\Lambda^{\le n}$ for the set
$\{\lambda \in \bigcup_{m \le n} \Lambda^m : d(\lambda) <
d(\lambda) + p \le n \implies r^{-1}(s(\lambda)) \cap \Lambda^p
= \emptyset\}$. Given a vertex $v$, a subset $F$ of $r^{-1}(v)$
is said to be exhaustive if for every $\mu \in r^{-1}(v)$ there
exists $\nu \in F$ such that $\MCE(\mu,\nu) \not= \emptyset$.
Given a finitely aligned $k$-graph $\Lambda$, a set
$\{s_\lambda : \lambda \in \Lambda\}$ of partial isometries is
called a Toeplitz-Cuntz-Krieger $\Lambda$-family if
\begin{itemize}
\item[(CK1)] $\{s_v : v \in \Lambda^0\}$ is a set of
    mutually orthogonal projections;
\item[(CK2)] $s_\mu s_\nu = s_{\mu\nu}$ whenever $s(\mu) =
    r(\nu)$; and
\item[(CK3)] $s^*_\mu s_\nu = \sum_{\mu\mu' = \nu\nu' \in
    \MCE(\mu,\nu)} s_{\mu'} s^*_{\nu'}$ for all $\mu,\nu
    \in \Lambda$.
\end{itemize}
It is called a Cuntz-Krieger $\Lambda$-family if it
additionally satisfies
\begin{itemize}
\item[(CK4)] $\prod_{\lambda \in F} (s_v - s_\lambda
    s^*_\lambda) = 0$ for all $v \in \Lambda^0$ and all
    nonempty finite exhaustive sets $F \subset r^{-1}(v)$.
\end{itemize}
The Cuntz-Krieger algebra $C^*(\Lambda)$ is the universal
$C^*$-algebra generated by a Cuntz-Krieger $\Lambda$-family.

Theorem~4.2 and Proposition~6.4 of \cite{RS1} show that Nica
covariant representations of $X(\Lambda)$ are in bijective
correspondence with Toeplitz-Cuntz-Krieger $\Lambda$-families.
To describe this bijection, we must briefly recall the
definition of $X(\Lambda)$ (see~\cite{RS1} for details). For $n
\in \NN^k$, we endow $c_c(\Lambda^n) = \lsp\{\delta_{\lambda} :
\lambda \in \Lambda^n\}$ with the structure of a pre-Hilbert
$c_0(\Lambda^0)$-bimodule via the following formulae: $(a \cdot
x \cdot b)(\lambda) := a(r(\lambda)) x(\lambda) b(s(\lambda))$;
and $\langle x, y \rangle_{c_0(\Lambda^0)}(v) =
\sum_{s(\lambda) = v} \overline{x(\lambda)}y(\lambda)$. Then
$X(\Lambda)_n$ is the completion of $c_c(\Lambda^n)$ in the
norm arising from $\langle
\cdot,\cdot\rangle_{c_0(\Lambda^0)}$. The isomorphisms
$X(\Lambda)_m \otimes X(\Lambda)_n \cong X(\Lambda)_{m+n}$ are
given by $\delta_\mu \otimes_{c_0(\Lambda^0)} \delta_\nu
\mapsto \delta_{\mu\nu}$ if $s(\mu) = r(\nu)$ (if $s(\mu) \not=
r(\nu)$, then $\delta_\mu \otimes_{c_0(\Lambda^0)} \delta_\nu =
0$ by definition of the balanced tensor product). Given a Nica
covariant representation $\psi$ of $X(\Lambda)$, the
corresponding Toeplitz-Cuntz-Krieger $\Lambda$-family (see
\cite[Definition~7.1]{RS1}) $\{t_\lambda : \lambda \in
\Lambda\}$ is defined by $t_\lambda = \psi(\delta_\lambda)$;
and we can recover $\psi$ from the $t_\lambda$ by linearity and
continuity.

\begin{prop}\label{prp:k-graph system}
Let $\Lambda$ be a finitely aligned $k$-graph, and let $X =
X(\Lambda)$ be the associated Cuntz-Krieger product system of
Hilbert bimodules. Let $\{s_\lambda : \lambda \in \Lambda\}$ be
the universal Cuntz-Krieger family in $C^*(\Lambda)$, and let
$j_X$ be the universal CNP-covariant representation of $X$ in
$\NO{X}$.
\begin{itemize}
\item[(1)] There is an isomorphism $\theta : C^*(\Lambda)
    \to \NO{X}$ satisfying $\theta(s_\lambda) =
    j_X(\delta_\lambda)$ for all $\lambda \in \Lambda$.
\item[(2)] Let $\psi : X \to B$ be a Nica-covariant
    representation of $X$, and let $\{t_\lambda : \lambda
    \in \Lambda\}$ be the corresponding
    Toeplitz-Cuntz-Krieger $\Lambda$-family $t_\lambda =
    \psi(\delta_\lambda)$. Then $\psi$ satisfies
    Definition~\ref{dfn:CP covariance} if and only if
    $\{t_\lambda : \lambda \in \Lambda\}$ satisfies the
    Cuntz-Krieger relation \cite[Definition~2.5(iv)]{RSY2}.
\end{itemize}
\end{prop}
\begin{proof}
Statement~(2) follows from~(1) and the universal properties of
$C^*(\Lambda)$ and $\NO{X}$. So it suffices to prove~(1).

For each $\lambda \in \Lambda$, let $j_\lambda =
j_X(\delta_\lambda)$, so $\{j_\lambda : \lambda \in \Lambda\}$
is a Toeplitz-Cuntz-Krieger $\Lambda$ family which generates
$\NO{X}$. We claim that $\{j_\lambda : \lambda \in \Lambda\}$
is a Cuntz-Krieger $\Lambda$-family.

Fix $v \in \Lambda^0$ and a finite exhaustive set $F \subset
v\Lambda$. We must show that
\[
\prod_{\mu \in F} (j_v - j_\mu j^*_\mu) = 0.
\]
For a subset $G \subset F$, we will denote by $\vee d(G)$ the
element $\bigvee_{\mu \in G} d(\mu)$ of $\NN^k$. If $\lambda,
\mu \in \Lambda$ satisfy $\lambda = \mu\mu'$ for some $\mu' \in
\Lambda$, we say that $\lambda$ extends $\mu$. Recall from
\cite{RS1} that for a nonempty subset $G$ of $F$, $\MCE(G)$
denotes the set $\{\lambda \in \Lambda : d(\lambda) = \vee
d(G), \lambda \text{ extends } \mu\text{ for all }\mu \in G\}$.
Recall also that $\vee F := \bigcup_{G \subset F} \MCE(G)$ is
finite and is closed under minimal common extensions. We have
\begin{align*}
\prod_{\mu \in F} (j_v - j_\mu j^*_\mu)
    &= j_v + \sum_{\substack{\emptyset \not= G \subset F \\ \lambda \in \MCE(G)}} (-1)^{|G|} j_\lambda j^*_\lambda \\
    &= j_X^{(e)}(\delta_v \otimes \delta_v^*)
        + \sum_{\substack{\emptyset \not= G \subset F \\ \lambda \in \MCE(G)}}
               (-1)^{|G|} j_X^{(\vee d(G))}(\delta_\lambda  \otimes \delta^*_\lambda).
\end{align*}
Since $j_{X(\Lambda)}$ satisfies~(CP) it suffices to show that
for each $q \in \NN^k$ there exists $r \ge q$ such that for all
$s \ge r$ we have
\[
\iotale^s_e(\delta_v \otimes \delta_v^*)
        + \sum_{\substack{\emptyset \not= G \subset F \\ \lambda \in \MCE(G)}}
               (-1)^{|G|} \iotale^s_{\vee d(G)}(\delta_\lambda \otimes \delta^*_\lambda) = 0.
\]
For this, fix $q \in \NN^k$, let $r = q \vee (\vee d(F))$ and
fix $s \ge r$. Since $\Xle{s} = \oplus_{t \le s}
\clsp\{\delta_\tau : \tau \in \Lambda^t \cap \Lambda^{\le s}\}$
(see Example~\ref{ex:Xle<->Lambdale}), it suffices to show that
for $\tau \in \Lambda^{\le s}$,
\begin{equation}\label{eq:at tau}
\Big(\iotale^s_e(\delta_v \otimes \delta_v^*)
        + \sum_{\substack{\emptyset \not= G \subset F \\ \lambda \in \MCE(G)}}
               (-1)^{|G|} \iotale^s_{\vee d(G)}(\delta_\lambda \otimes \delta^*_\lambda)\Big)(\delta_\tau) = 0.
\end{equation}
Fix $\tau \in \Lambda^{\le s}$. For any $\mu \in F$, we have $s \ge
d(\mu)$, so
\begin{equation}\label{eq:action on tXs}
\iotale^s_{d(\mu)}(\delta_\mu \otimes \delta^*_\mu)(\delta_\tau) =
\begin{cases}
\delta_\tau &\text{ if $\tau$ extends $\mu$}\\
0 &\text{ otherwise}.
\end{cases}
\end{equation}

Fix a nonempty subset $G$ of $F$. Then
\[
\Big(\prod_{\mu \in G} \iotale^s_{d(\mu)}(\delta_\mu \otimes \delta^*_\mu)\Big)(\delta_\tau) =
\begin{cases}
\delta_\tau &\text{ if $\tau$ extends each $\mu$ in $G$}\\
0 &\text{ otherwise}.
\end{cases}
\]
The factorisation property implies that $\tau$ extends each
$\mu$ in $G$ if and only if there exists $\lambda$ in $\MCE(G)$
such that $\tau$ extends $\lambda$. The factorisation property
also implies that if there does exist such a $\lambda \in
\MCE(G)$ then it is necessarily unique. We therefore have
\[
\Big(\prod_{\mu \in G} \iotale^s_{d(\mu)}(\delta_\mu \otimes \delta^*_\mu)\Big)(\delta_\tau) =
\Big(\sum_{\lambda \in \MCE(G)} \iotale^s_{\vee d(G)}(\delta_\lambda \otimes \delta^*_\lambda)\Big)(\delta_\tau)
\]

Since the fixed nonempty subset $G$ of $F$ in the preceding
paragraph was arbitrary, we may now calculate:
\begin{align*}
\Big(\prod_{\mu \in F} \big(\iotale^s_e(\delta_v \otimes \delta^*_v)
 &- \iotale^s_{d(\mu)} (\delta_\mu \otimes \delta^*_\mu)\big)\Big)(\delta_\tau) \\
 &=  \Big(\iotale^s_e(\delta_v \otimes \delta_v^*)
        + \sum_{\emptyset \not= G \subset F}
               \Big((-1)^{|G|} \prod_{\mu \in G} \iotale^s_{d(\mu)}(\delta_\mu \otimes \delta^*_\mu)\Big)\Big)(\delta_\tau) \\
 &= \Big(\iotale^s_e(\delta_v \otimes \delta_v^*)
        + \sum_{\substack{\emptyset \not= G \subset F \\ \lambda \in \MCE(G)}}
               (-1)^{|G|} \iotale^s_{\vee d(G)}(\delta_\lambda  \otimes \delta^*_\lambda)\Big)(\delta_\tau).
\end{align*}
Since $F$ is exhaustive, there exists $\nu \in F$ such that
$\MCE(\tau, \nu) \not=\emptyset$; say $\tau\tau' \in
\MCE(\tau,\nu$). Since $d(\nu), d(\tau) \le s$, we have
$d(\tau\tau') \le s$. Since $\tau \in \Lambda^{\le s}$, this
forces $d(\tau') = 0$, so $\tau$ extends $\nu$.
Thus~\eqref{eq:action on tXs} implies that
\begin{align*}
\Big(\prod_{\mu \in F}& \big(\iotale^s_e(\delta_v \otimes \delta^*_v) -
    \iotale^s_{d(\mu)} (\delta_\mu \otimes \delta^*_\mu)\big)\Big)(\delta_\tau) \\
    &= \Big(\prod_{\mu \in F\setminus\{\nu\}} (\iotale^s_e(\delta_v \otimes \delta^*_v) -
    \iotale^s_{d(\mu)} (\delta_\mu \otimes \delta^*_\mu))\Big) \big(
    \big(\iotale^s_e(\delta_v \otimes \delta^*_v) -
    \iotale^s_{d(\nu)} (\delta_\nu \otimes \delta^*_\nu)\big)(\delta_\tau)
    = 0,
\end{align*}
establishing~\eqref{eq:at tau}. Hence $\{j_\lambda : \lambda \in
\Lambda\}$ is a Cuntz-Krieger $\Lambda$-family as claimed.

Since the $j_\lambda$ generate $\NO{X}$, the universal property
of $C^*(\Lambda)$ implies that there is a surjective
homomorphism $\theta : C^*(\Lambda) \to \NO{X}$ satisfying
$\theta(s_\lambda) = j_\lambda = j_X(\delta_\lambda)$ for all
$\lambda \in \Lambda$. By Theorem~\ref{thm:isometric}, we have
$j_v \not= 0$ for all $v \in \Lambda^0$. Since $\theta$
intertwines the gauge actions of $\TT^k$ on $\NO{X}$ and
$C^*(\Lambda)$, the gauge-invariant uniqueness theorem for
$C^*(\Lambda)$ \cite[Theorem~4.2]{RSY2} therefore implies that
$\theta$ is an isomorphism.
\end{proof}

\subsection{Boundary quotients of Toeplitz algebras}\label{sec:CrispLaca}

In this section, we consider product systems whose fibres are
isomorphic to ${_\CC}\CC_\CC$. Nica covariant representations
of such product systems amount to Nica covariant
representations of $(G,P)$ in the sense of \cite{CL1, LR,
Nica}. This prompts us to explore the connection between our
$\NO{X}$ and the boundary quotients of $\Tt(G,P)$ studied by
Crisp and Laca in \cite{CL}.

The first part of the following proposition follows from
results of Fowler and Raeburn \cite{FR98}, but we include it
for completeness. We first make the following definition: if
$(G,P)$ is a quasi-lattice ordered group, we say that a finite
subset $F$ of $P$ is a \emph{foundation set} for $P$ if, for
every $q \in P$ there exists $p \in F$ such that $p \vee q <
\infty$.

Note that what we have called foundation sets are precisely the
\emph{boundary relations} of \cite[Definition~3.4]{CL}.

\begin{notation}\label{ntn:C^P}
Given a semigroup $P$, we denote by $\CC^P$ the product system
of right-Hilbert \mbox{$\CC$--$\CC$} bimodules determined by
$X_p = {_\CC\CC_\CC}$ for all $p$ with multiplication in $X$
given by multiplication of complex numbers. As a notational
convenience, when we are regarding a complex number $z$ as an
element of $\CC^P_p$, we shall denote it $z_p$.
\end{notation}

\begin{prop}\label{prp:semigroup universal property}
Let $(G,P)$ be a quasi-lattice ordered group. Let $X = \CC^P$
be the product system discussed above. Then $\Tc(X)$ is the
universal $C^*$-algebra generated by a semigroup representation
$p \mapsto V_p$ of $P$ by isometries satisfying
\begin{equation}\label{eq:semigroup Nica}
V_p V^*_p V_q V^*_q =
\begin{cases}
V_{p \vee q} V^*_{p \vee q} &\text{ if $p \vee q < \infty$} \\
0 &\text{ otherwise.}
\end{cases}
\end{equation}
The images $\{W_p : p \in P\}$ of these isometries under the
canonical homomorphism $(j_X)_* : \Tc(X) \to \NO{X}$ satisfy
the additional relation
\begin{equation}\label{eq:semigroup CP}
\prod_{p \in F} (1 - W_p W^*_p) = 0\text{ for every foundation set $F$ for $P$.}
\end{equation}
\end{prop}
\begin{proof}
It is straightforward to check that the elements $V_p =
i_X(1_p)$ of $\Tc(X)$ generate $\Tc(X)$, determine a semigroup
representation of $P$, and satisfy~\eqref{eq:semigroup Nica}.
We therefore need only show that their images $W_p = j_X(1_p)$
in $\NO{X}$ satisfy~\eqref{eq:semigroup CP}.

For this, fix a foundation set $F$ for $P$. In what follows, we
shall write $\vee H$ for the least upper bound of a finite
subset $H$ of $P$ when it exists, and when it does not exist,
we shall write $\vee H = \infty$. Since the $W_p$ are the
images of the $V_p$ under a homomorphism,
relation~\eqref{eq:semigroup Nica} holds amongst the $W_p$. We
have
\begin{align*}
\prod_{p \in F} (1 - W_p W^*_p)
    &= 1 + \sum_{\substack{\emptyset \not= H \subset F \\ \vee H < \infty}} (-1)^{|H|} W_{\vee H} W^*_{\vee H} \\
    &= 1 + \sum_{\substack{\emptyset \not= H \subset F \\ \vee H < \infty}} (-1)^{|H|} j_X^{(\vee H)}(1_{\vee H} \otimes 1_{\vee H}^*).
\end{align*}
Since $j^{(e)}(1_e \otimes 1^*_e) = W_e = 1$, since
$\iotale^s_e(1_e \otimes 1^*_e) = 1_{\Ll(\Xle{s})}$, and since
$j_X$ is Cuntz-Pimsner covariant, it suffices to show that
\[
1_{\Ll(\Xle{s})} + \sum_{\substack{\emptyset \not= H \subset F \\ \vee H < \infty}} (-1)^{|H|}
\iotale_{\vee H}^s(1_{\vee H} \otimes 1_{\vee H}^*) = 0\text{ for large $s$.}
\]
For this, fix $q \in P$. Since $F$ is a foundation set for $P$,
we must have $q \vee a < \infty$ for some $a \in F$. An
inductive argument then shows that there exists $r \in P$ such
that $r \ge q \vee a$ and such that for each $p \in F$ either
$r \ge p$ or $r \vee p = \infty$. Fix $s \in P$ with $s \ge r$.

Fix $\emptyset \not= H \subset F$ with $\vee H < \infty$. Then
\[
\iota^s_{\vee H}(1_{\vee H} \otimes 1^*_{\vee H}) =
\begin{cases}
1_{\Ll(X_s)} &\text{ if $\vee H \le s$}\\
0 &\text{ otherwise,}
\end{cases}
\]
and for $p \in H$, we have
\[
\iota^s_p(1_p \otimes 1^*_p) =
\begin{cases}
1_{\Ll(X_s)} &\text{ if $p \le s$}\\
0 &\text{ otherwise.}
\end{cases}
\]
Since $\vee H \le s$ if and only if $p \le s$ for all $p \in H$, the
preceding two displayed equations imply that
\[
\iota^s_{\vee H}(1_{\vee H} \otimes 1_{\vee H}^*)
 = \prod_{p \in H} \iota^s_p(1_p \otimes 1^*_p).
\]

Since each $\phi_t$ is injective, $\Xle{p} = X_p$ and
$\iotale^s_p(1_p \otimes 1^*_p) = \iota^s_p(1_p \otimes 1^*_p)$
for all $p,s \in P$. Since the calculations in the previous
paragraph are valid for arbitrary nonempty $H \subset F$ with
$\vee H < \infty$, we may now calculate
\[\begin{split}
1_{\Ll(\Xle{s})} + \sum_{\substack{\emptyset \not= H \subset F \\ \vee H \not= \infty}}
    (-1)^{|H|} \iotale_{\vee H}^s(1_{\vee H} \otimes 1_{\vee H}^*)
 &= 1_{\Ll(X_s)} + \sum_{\substack{\emptyset \not= H \subset F \\ \vee H \not= \infty}}
    (-1)^{|H|} \iota_{\vee H}^s(1_{\vee H} \otimes 1_{\vee H}^*) \\
 &= \prod_{p \in F} (1_{\Ll(X_s)} - \iota^s_p(1_p \otimes 1^*_p)) \\
 &= \begin{cases}
    1_{\Ll(X_s)} &\text{ if $p \not\le s$ for all $p \in F$}\\
    0 &\text{ otherwise.}
 \end{cases}
\end{split}\]
By choice of $r$ and $s$, we have $a \le s$ for some $a \in F$,
so
\[
1_{\Ll(\Xle{s})} + \sum_{\substack{\emptyset \not= H \subset F \\ \vee H \not= \infty}}
    (-1)^{|H|} \iotale_{\vee H}^s(1_{\vee H} \otimes 1_{\vee H}^*) = 0
\]
as required.
\end{proof}

The point of the above proposition is the relationship it
suggests with Crisp and Laca's boundary quotient of $\Tt(G,P)$
(see \cite[Definitions 3.1~and~3.4]{CL}).

Recall that given an undirected graph $\Gamma$ with vertex set
$S$, the right-angled Artin group $G$ associated to $\Gamma$ is
the group $G = \langle S \,|\, st = ts$ whenever $s$ and $t$
are adjacent in $\Gamma\rangle$. We write $P$ for the submonoid
of $G$ generated by $S$. The pair $(G,P)$ is a quasi-lattice
ordered group and satisfies~\eqref{eq:max elements}.

For this $(G,P)$, \cite[Theorem~6.7]{CL} shows that the
boundary quotient $C_0(\partial\Omega) \rtimes G$ is simple and
is the universal $C^*$-algebra generated by elements $\{T_s : s
\in S\}$ satisfying
\begin{itemize}
\item[(1)] $T_s^* T_s = 1$ for all $s \in S$;
\item[(2)] $T_s T_t = T_t T_s$ and $T^*_s T_t = T_t T^*_s$
    whenever $s$ and $t$ are adjacent in $\Gamma$;
\item[(3)] $T^*_s T_t = 0$ whenever $s$ and $t$ are
    adjacent in $\Gopp$; and
\item[(4)] $\prod_{s \in C} (1 - T_s T^*_s) = 0$ when $C$
    is the vertex set of any finite connected component of
    $\Gopp$.
\end{itemize}

\begin{cor}
Let $\Gamma$ be an undirected graph with vertex set $S$, and
let $(G,P)$ be the associated right-angled Artin group. Suppose
that $G$ has trivial centre. Let $X = \CC^P$ be the product
system of Notation~\ref{ntn:C^P}. Let $C_0(\partial\Omega)
\rtimes G$ be the boundary quotient of~\cite{CL}. Then $\NO{X}
\cong C_0(\partial\Omega) \rtimes G$.
\end{cor}
\begin{proof}
As in \cite{CL}, let $\Gopp$ denote the graph which has the
same vertex set $S$ such that $s$ and $t$ are adjacent in
$\Gopp$ if and only if they are not adjacent in $\Gamma$.

For $p \in P$, let $W_p = j_X(1_p) \in \NO{X}$.
Proposition~\ref{prp:semigroup universal property} shows that
$p \mapsto W_p$ is a semigroup representation of $P$ by
isometries which satisfy \eqref{eq:semigroup
Nica}~and~\eqref{eq:semigroup CP}. Since $S$ generates $P$, the
set $\{W_s : s \in S\}$ generates $\NO{X}$ and satisfies~(1)
and the first part of~(2). If $s,t$ are adjacent in $\Gamma$,
we have $s \vee t = st = ts$, and so \eqref{eq:semigroup Nica},
(1) and the first part of~(2) force $W^*_s W_t = W^*_s (W_s
W^*_s W_t W^*_t) W_t = W^*_s W_{st} W^*_{ts} W_t = W_t W^*_s$,
so the $W_s$ satisfy~(2). If $s,t$ are adjacent in $\Gopp$,
then~\eqref{eq:semigroup Nica} gives $W^*_s W_t = W^*_s (W_s
W^*_s W_t W^*_t) W_t = 0$. Hence the $W_s$ satisfy (1)--(3).
The final paragraph of the proof of \cite[Theorem~6.7]{CL} and
\cite[Definition~3.4]{CL} imply that the vertex set $C$ of any
finite connected component of $\Gopp$ is a foundation set for
$P$. Hence Proposition~\ref{prp:semigroup universal property}
implies that the $W_s$ satisfy~(4).

The universal property of $C(\partial\Omega) \rtimes G$ now
implies that there is a homomorphism $\pi : C(\partial\Omega)
\rtimes G \to \NO{X}$ satisfying $\pi(T_s) = W_s$ for all $s
\in S$. Since $S$ generates $P$ and since $\{W_p : p \in P\}$
generates $\NO{X}$, $\pi$ is surjective.
Theorem~\ref{thm:isometric} implies that $\NO{X} \not= \{0\}$.
Since $C_0(\partial\Omega) \rtimes G$ is simple, it follows
that $\pi$ is an isomorphism.
\end{proof}


\begin{thebibliography}{00}

\bibitem{Arv} W. Arveson, \emph{Continuous analogues of Fock
    space}, Memoirs Amer. Math. Soc. {\bf 80} (1989).

\bibitem{BHRS} T. Bates, J. Hong, I. Raeburn, and W. Szyma\'nski,
    \emph{The ideal structure of the $C^*$-algebras of infinite
    graphs}, Illinois J. Math. {\bf 46} (2002), 1159--1176.

\bibitem{Black} B. Blackadar, \emph{Operator algebras. Theory of
    $C\sp *$-algebras and von Neumann algebras}. Encyclopaedia of
    Mathematical Sciences, 122. Operator Algebras and Non-commutative
    Geometry, III. Springer-Verlag, Berlin, 2006.

\bibitem{CLSV} T. Carlsen, N. Larsen, A. Sims and S.
    Vittadello, \emph{A gauge-invariant uniqueness theorem for Cuntz-Pimsner
    algebras of product systems of Hilbert bimodules}, in
    preparation, 2007.

\bibitem{Coburn:isometryI}
    L.A. Coburn, \emph{The $C^*$-algebra generated by an isometry.}
    Bull. Amer. Math. Soc. {\bf73} (1967), 722--726.

\bibitem{CL1} J. Crisp and M. Laca, \emph{On the Toeplitz algebras of
    right-angled and finite-type Artin groups}, J. Austral. Math.
    Soc. {\bf72} (2002), 223--245.

\bibitem{CL} J. Crisp and M. Laca, \emph{Boundary quotients
    and ideals of Toeplitz $C^*$-algebras of Artin groups}, J.
    Funct. Anal. {\bf242} (2007), 127--156.

\bibitem{CK} J. Cuntz and W. Krieger, \emph{A class of $C^*$-algebras
    and topological Markov chains}, Invent.  Math. {\bf 56} (1980),
    251--268.

\bibitem{Din} H. T. Dinh, \emph{Discrete product systems and
    their
    $C^*$-algebras}, J. Funct. Anal. {\bf 102} (1991), 1--34.

\bibitem{Din2} H. T. Dinh, \emph{On generalized Cuntz
    $C^*$-algebras}, J. Operator Theory {\bf 30} (1993), 123--135.

\bibitem{EL} R. Exel and M. Laca, \emph{Cuntz-Krieger algebras
    for infinite matrices}, J. reine angew. Math. {\bf 512} (1999),
    119--172.

\bibitem{F99} N. J. Fowler, \emph{Discrete product systems of Hilbert
    bimodules}, Pacific J. Math.  {\bf 204} (2002), 335--375.

\bibitem{FLR} N. J. Fowler, M. Laca, and I. Raeburn, \emph{The
    $C^*$-algebras of infinite graphs}, Proc.  Amer.  Math. Soc.
    {\bf 128} (2000), 2319--2327.

\bibitem{FMR} N. J. Fowler, P. S. Muhly, and I. Raeburn,
    \emph{Representations of Cuntz-Pimsner algebras}, Indiana Univ.
    Math. J. {\bf 52} (2003), 569--605.

\bibitem{FR98} N. J. Fowler and I. Raeburn, \emph{Discrete product
    systems and twisted crossed products by semigroups}, J. Funct.
    Anal. {\bf 155} (1998), 171--204.

\bibitem{FR} N. J. Fowler and I. Raeburn, \emph{The Toeplitz
    algebra of a Hilbert bimodule}, Indiana Univ.  Math.  J. {\bf 48} (1999),
    155--181.

\bibitem{K5} T. Katsura, \emph{A construction of
    $C^*$-algebras from
    $C^*$-correspondences}, Advances in Quantum Dynamics, 173--182,
    Contemp. Math, {\bf335}, Amer. Math. Soc., Providence, RI, 2003.

\bibitem{K6} T. Katsura, \emph{On $C^*$-algebras associated
    with
    $C^*$-correspondences}, J. Funct. Anal. {\bf217} (2004),
    366--401.

\bibitem{KP} A. Kumjian and D. Pask, \emph{Higher rank graph
    $C^*$-algebras}, New York J. Math.  {\bf 6} (2000), 1--20.

\bibitem{KPR} A. Kumjian, D. Pask, and I. Raeburn,
    \emph{Cuntz-Krieger algebras of directed graphs}, Pacific J.
    Math. {\bf 184} (1998), 161--174.

\bibitem{KPRR} A. Kumjian, D. Pask, I. Raeburn, and J.
    Renault,
    \emph{Graphs, groupoids and Cuntz-Krieger algebras}, J. Funct.
    Anal.  {\bf 144} (1997), 505--541.

\bibitem{LR} M. Laca and I. Raeburn, \emph{Semigroup crossed
    products and the Toeplitz algebras of nonabelian groups},
    J. Funct. Anal. {\bf139} (1996), 415--440.

\bibitem{Lan} E. C. Lance, \emph{Hilbert $C^*$-modules: A
    toolkit for operator algebraists}, London Math. Soc. Lecture Note Series,
    vol. 210, Cambridge Univ. Press, Cambridge, 1994.

\bibitem{MS} P.S. Muhly and B. Solel, \emph{Tensor algebras
    over $C^*$-correspondences (representations, dilations, and
    $C^*$-envelopes)}, J. Funct. Anal. {\bf 158} (1998), 389--457.

\bibitem{Nica} A. Nica, \emph{$C^*$-algebras generated by
    isometries
    and Wiener-Hopf operators}, J. Operator Theory, {\bf 27} (1992),
    17--52.

\bibitem{Pim97} M. V. Pimsner, \emph{A class of $C^*$-algebras
    generalizing both Cuntz-Krieger algebras and crossed products by
    $\mathbb{Z}$}, Fields Institute Communications {\bf 12} (1997),
    189--212.

\bibitem{CBMSbk} I. Raeburn, \emph{Graph algebras}, CBMS
    Regional Conference Series in Mathematics, Vol. 103, Amer. Math. Soc.,
    2005.

\bibitem{RS1} I. Raeburn and A. Sims, \emph{Product systems of
    graphs
    and the Toeplitz algebras of higher-rank graphs}, J. Operator Th.
    {\bf53} (2005), 399--429.

\bibitem{RSY1} I. Raeburn, A. Sims, and T. Yeend,
    \emph{Higher-rank
    graphs and their $C^*$-algebras}, Proc.  Edinb. Math.  Soc.  {\bf
    46} (2003), 99--115.

\bibitem{RSY2} I. Raeburn, A. Sims, and T. Yeend, \emph{The
    $C^*$-algebras of finitely aligned higher-rank graphs}, J. Funct.
    Anal. {\bf213} (2004), 206--240.

\bibitem{TFB} I. Raeburn and D. P. Williams, \emph{Morita
    equivalence
    and continuous-trace $C^*$-algebras}, Math. Surveys and
    Monographs, vol. 60, Amer. Math. Soc., Providence, 1998.

\end{thebibliography}
\end{document}